%% file: main.tex
\documentclass[11pt]{article}
    \pdfoutput=1
	
	\newcommand{\blind}{0}
	
    \makeatletter
    \renewcommand\section{\@startsection {section}{1}{\z@}%
                                       {-3.5ex \@plus -1ex \@minus -.2ex}%
                                       {2.3ex \@plus.2ex}%
                                       {\normalfont\fontfamily{phv}\fontsize{16}{19}\bfseries}}
    \renewcommand\subsection{\@startsection{subsection}{2}{\z@}%
                                          {-3.25ex\@plus -1ex \@minus -.2ex}%
                                         {1.5ex \@plus .2ex}%
                                         {\normalfont\fontfamily{phv}\fontsize{14}{17}\bfseries}}
    \renewcommand\subsubsection{\@startsection{subsubsection}{3}{\z@}%
                                        {-3.25ex\@plus -1ex \@minus -.2ex}%
                                         {1.5ex \@plus .2ex}%
                                         {\normalfont\normalsize\fontfamily{phv}\fontsize{14}{17}\selectfont}}
    \makeatother
	
	\usepackage{amsmath}
	\usepackage{graphicx}
	\usepackage{enumerate}
	\usepackage{xcolor}
	\usepackage{natbib} 
	\usepackage{url} 
	
    \usepackage[utf8]{inputenc}
    \usepackage{algorithm}
    \usepackage{algpseudocode}
    \usepackage{amssymb}
    \usepackage{accents}
    \usepackage{array}
    \usepackage{enumitem}
    \usepackage[T1]{fontenc}
    \usepackage[margin=1.00in]{geometry}
    \usepackage{layouts}  
    \usepackage{lipsum}
    \usepackage{multirow,multicol}
    \usepackage[short,nocomma]{optidef}
    \usepackage{subcaption}
    \usepackage{hyperref}

\DeclareMathOperator{\prob}{Pr}
\def\tigerdam{Tiger Dam\texttrademark}
\def\tigerdams{Tiger Dams\texttrademark}
\def\wt{\widetilde}
\def\bs{\boldsymbol}
\newcolumntype{L}[1]{>{\raggedright\let\newline\\\arraybackslash\hspace{0pt}}m{#1}}
\newcolumntype{C}[1]{>{\centering\let\newline\\\arraybackslash\hspace{0pt}}m{#1}}
\newcolumntype{R}[1]{>{\raggedleft\let\newline\\\arraybackslash\hspace{0pt}}m{#1}}
\let\emptyset\varnothing

\DeclareMathOperator*{\argmax}{argmax}

\let\epsswitch\undefined  

\definecolor{response}{RGB}{0, 0, 0}
 
\begin{document}
\def\spacingset#1{\renewcommand{\baselinestretch}%
                               {#1}\small\normalsize} \spacingset{1}
\if0\blind
{
    \title{\bf A two-stage stochastic programming model for electric substation flood mitigation prior to an imminent hurricane}
    \author{Brent Austgen$^*$, Erhan Kutanoglu, and John J. Hasenbein \\
            Operations Research and Industrial Engineering Program, \\
            The University of Texas at Austin, Austin, TX, USA \\
            *Corresponding Author: brent.austgen@utexas.edu}
    \date{}
    \maketitle
}
\fi
\if1\blind
{
    \title{\bf \emph{IISE Transactions} \LaTeX \ Template}
    \author{Author information is purposely removed for double-blind review}
    \bigskip
    \bigskip
    \bigskip
    \begin{center}
        {\LARGE\bf A two-stage stochastic programming model for electric substation flood mitigation prior to an imminent hurricane}
    \end{center}
    \medskip
}
\fi
\bigskip
\begin{abstract}
We present a stochastic programming model for informing the deployment of ad hoc flood mitigation measures to protect electrical substations prior to an imminent and uncertain hurricane. The first stage captures the deployment of a fixed number of mitigation resources, and the second stage captures grid operation in response to a contingency. The primary objective is to minimize expected load shed. We develop methods for simulating flooding induced by extreme rainfall and construct two geographically realistic case studies, one based on Tropical Storm Imelda and the other on Hurricane Harvey. Applying our model to those case studies, we investigate the effect of the mitigation budget on the optimal objective value and solutions. Our results highlight the sensitivity of the optimal mitigation to the budget, a consequence of those decisions being discrete. We additionally assess the value of having better mitigation options and the spatial features of the optimal mitigation.
\end{abstract}
\noindent%
{\it Keywords:} stochastic programming, discrete optimization, power grid, flooding, mitigation, resilience, risk management
\spacingset{1.5} 

\input{introduction/_main}
\input{review/_main}
\input{modeling/_main}
\input{casestudies/_main}
\input{results/_main}
\input{conclusions/_main}

\if0\blind{
\section*{Acknowledgements}
The authors thank The University of Texas Energy Institute for funding our ``Defending the Electricity Infrastructure Against Extreme Weather Events, Now \& in the Future'' research project as part of the ``Fueling a Sustainable Energy Transition'' initiative.

The authors also acknowledge the Texas Advanced Computing Center (TACC) at The University of Texas at Austin for providing HPC and data visualization resources that have contributed to the research results reported within this paper. URL: http://www.tacc.utexas.edu
\fi

\bibliographystyle{chicago}
\spacingset{1}
\bibliography{main}

\end{document}

%% file: introduction/_main.tex
\section{Introduction} \label{section:intro}

Since 1980, 338 natural disasters each costing at least \$1B have occurred in the United States \citep{BillionDollarDisasters2022}.
Though only 59 of these 338 disasters were tropical cyclone (TC) events, the 7 most costly disasters were all TCs as were 26 of the 50 most costly disasters.
The annual cost of high-profile TCs has generally increased with 2005 (Hurricane Katrina) and 2017 (Hurricanes Harvey, Irma, and Maria) leading all other years by a sizable margin.
Many climate models project the frequency and intensity of the most extreme TCs (\textit{i.e.}, Category 4 and 5 storms) to increase \citep{Webster2005,Knutson2020}. Granted, those same models project a decrease in the total number of TCs globally except in the North Atlantic basin.

Natural disasters affect lifeline infrastructure like water, natural gas, transportation, and communication systems, emergency and healthcare services, gas stations, and grocery stores. These systems are all either dependent on or codependent with the power grid. Power grid resilience is thus key to the resilience of modern infrastructure; however, extreme weather has repeatedly proven to be a major threat to the power grid. Hurricane Harvey, for example, was a particularly devastating storm. According to the North American Electric Reliability Corporation (NERC), the storm affected more than 200 transmission lines and over 200 load-serving substations in the domain of the Electric Reliability Council of Texas (ERCOT) alone \citep{NERC2018}. In total, over 2 million customers were affected.

\textcolor{response}{In NERC's report, the potential for substation flooding near the Texas coast was identified as a concern, but no preparation activities were performed to mitigate that flooding \citep{NERC2018}. To understand why this is a problem, consider that the power grid operates at two levels: transmission and distribution. The transmission grid is responsible for transporting high-voltage electricity over long distances from power plants to substations. The distribution grid operates at lower voltages, distributing electricity from substations to local customers. Substations serve as the interface between the transmission grid and the distribution grid. Specifically, they facilitate the flow of power throughout the transmission grid and step down the electrical voltage to levels suitable for satisfying demand on the distribution grid. When flooding affects a substation, critical components such as transformers or control houses may be damaged \citep{Boggess2014}, resulting in prolonged power outages for customers served by the affected substation, potentially including tens of thousands of residents, and causing atypical grid operation that impacts customers elsewhere. To help prevent outages, we consider in this paper the deployment of a fixed number of \tigerdam~flood mitigation resources to protect vulnerable substations.}

\textcolor{response}{To inform mitigation decisions prior to an imminent but still uncertain hurricane, we propose using a two-stage stochastic programming model. In our model, the first-stage decisions pertain to the flood mitigation decisions. For each flooding scenario considered, the associated second-stage problem is to operate the potentially degraded transmission grid in the wake of the disaster with the primary objective of minimizing unserved load (\textit{i.e.}, customer outages). Transmission grid operations may be represented by the alternating current (AC) power flow model; however, the model is nonconvex and thus computationally prohibitive in stochastic power system planning models. Accordingly, we model transmission grid operations in the second-stage problem using the widely used direct current (DC) approximation, a linear surrogate of ACPF. For further details on the modeling of power grid operations, including AC and DC power flow models, we refer the reader to \citet{Frank2016}, \citet{Molzahn2019}, and \citet{Bienstock2015}.}

In this paper, we deliver the following contributions.
\begin{enumerate}
    \item We develop a two-stage stochastic programming model for informing flood mitigation decision making prior to an imminent and uncertain hurricane. In the first stage, we incorporate a novel model for the resource- and technology-constrained deployment of temporary flood barriers. Our model of flood mitigation improves upon existing models \citep{Souto2022,Movahednia2022b} by accounting for the ability to stack flood barriers and the potential for deployed mitigation to be insufficient. In the second stage, we leverage the DC power flow approximation of system operation to capture the consequences of both flood barrier deployment decisions and unmitigated substation flooding on the power grid.
    \item We introduce novel scenario generation methods that leverage the uncertain but still quantifiable aspects of a TC and a state-of-the-art stream flow model to generate flooding scenarios. We apply these methods to Hurricane Harvey and Tropical Storm Imelda data to construct two geographically realistic case studies that augment a widely used 2000-bus synthetic grid set on the footprint of Texas. Our efforts further enable the research of realistic weather impact on power grids without requiring real grid data.
    \item We present a parametric greedy heuristic for quickly obtaining multiple mitigation solutions, which we leverage to warmstart the optimization solver (Gurobi). \textcolor{response}{For every studied instance, we computed the optimal objective value. Using these values, we show that the best heuristic mitigation solution for each instance achieves an objective value within 5\% of the optimal solution.}
    \textcolor{response}{\item Applying our model to the two case studies,} we investigate the impact of the mitigation budget on the optimal objective value and mitigation solution \textcolor{response}{in a sensitivity study. Our results highlight a transitory phenomenon in the optimal mitigation solutions that arises from the discrete nature of the mitigation decision making. We additionally assess the effect of the technological limit of mitigation on the optimal mitigation solution with specific focus on siting and sizing.}
\end{enumerate}

The remainder of the paper is structured as follows. Section~\ref{section:review} provides a review of literature pertaining to power grid resilience decision making. In Section~\ref{section:modeling}, we define the two-stage model that we propose be used to inform flood mitigation decision making for the transmission grid. In Section~\ref{section:case_studies}, we discuss the development of the geographically realistic case studies. In Section~\ref{section:results}, we introduce our solution methodology and present the results of applying our model to the two case studies. Finally, we present our conclusions in Section~\ref{section:conclusions}.

%% file: review/_main.tex
\section{Literature Review} \label{section:review}

\subsection{\emph{Risk and Resilience Modeling}}

A variety of definitions for ``risk'' and ``resilience'' exist in the literature. One of the most popular and comprehensive definitions of resilience comes from a report by the National Research Council: ``resilience is the ability to prepare and plan for, absorb, recover from, and more successfully adapt to adverse events'' \citep{NRC2012}. More nuanced aspects of risk and resilience are debated. For example, risk and resilience as defined in \cite{Linkov2019} are threat-dependent and threat-agnostic, respectively. However, in \cite{Logan2022}, resilience is defined as being system- and context-specific and inherently integrated with risk such that the two are only reasonably assessed together. We adopt the latter perspective in our work and develop our model according to the conceptual framework of \cite{Watson2014} that proposes resilience assessments be formed with consideration of three key factors: the threat, the likelihood, and the consequences.

For resilience via mitigation, decisions typically are framed with respect to one of two broad classes of threats: designed or natural. Threats from the former class are characterized by malicious intent. For modeling mitigation, appropriate decision models for this type of threat include bilevel and multilevel optimization models \citep{GhorbaniRenani2020} and robust optimization models \citep{Hien2020,Borrero2021}. Threats from the latter class are characterized by randomness (\textit{i.e.}, a probability distribution). Suitable models for mitigation against natural threats include stochastic optimization models \citep{Zhang2023} and distributionally robust optimization models \citep{Belle2023}.

In this paper, we employ stochastic programming to model mitigation decisions for hurricanes, a natural threat. We capture the threat and likelihoods by considering a sample distribution of representative hurricane flooding scenarios. The consequences in each scenario are a function of the flooding realization, the mitigation decisions, and the response decisions.

\subsection{\emph{Power Flow Modeling}}

There are many examples of surrogate PF models being incorporated into optimization models for power grid resilience to natural disasters. Examples include the network flow relaxation for general distribution system hardening \citep{Tan2018} and grid operation during a progressing wildfire \citep{Mohagheghi2015}; the DC power transfer distribution factor (PTDF) approximation for winter storm mitigation planning \citep{Garcia2022}; the DC B-theta approximation for hurricane restoration planning \citep{Arab2015}, hurricane mitigation planning \citep{Shukla2022,Movahednia2022b}, winter storm mitigation planning \citep{Pierre2018}, inventory stockpiling \citep{Coffrin2011}, and proactive grid posturing \citep{SahraeiArdakani2017,Quarm2022}; the LPAC approximation for transmission system restoration \citep{Coffrin2015}; and the second-order cone programming (SOCP) relaxation for transmission system hardening \citep{Garifi2022}.

In the aforementioned research, most grid instances to which optimization models were applied comprised around 100 buses. The largest grid instances studied were the ACTIVS 2000-bus synthetic grid of Texas \citep{Birchfield2017} and a 1263-bus grid representative of Puerto Rico \citep{Elizondo2020} used by \citet{Garcia2022} and \citet{Quarm2022}, respectively. The models in these papers incorporate linear DC approximations of PF but have a static grid topology in each contingency scenario. In this paper, we appropriate an adaptation of the DC approximation that allows the grid topology to vary as a function of both the contingency and the resilience decisions as in \citet{Pierre2018}, \citet{Movahednia2022b}, and \citet{Garifi2022} and furthermore allows total blackouts to occur for the sake of model feasibility. We apply our PF model to the coast-focused reduction of the ACTIVS 2000-bus grid comprising 663 buses \citep{Shukla2022,Austgen2022b,Souto2022}.

\subsection{\emph{Hurricane Modeling}}

Much of the research at the intersection of hurricanes and power grids centers around wind-related damages. A common approach for generating wind contingencies is to use historical or projected hurricane wind speed data in tandem with component fragility curves. This approach is followed in \citet{Sabouhi2020} and \citet{Bennett2021} to evaluate power grid resilience and in \citet{Quarm2022} and \citet{Poudel2019} to inform resilience decision making for power grids. In \citet{Mensah2014}, wind fragility curves are incorporated in a Bayesian network model driven by PF modeling to predict outages caused by hurricane winds. Data-driven predictive approaches such as accelerated failure time models \citep{Liu2007} and random forest models \citep{Guikema2014} have also been used to assess wind-oriented hurricane consequences.

Flood-related power grid damages are less prevalent in the literature. In \citet{Shukla2022}, damages induced by storm surge are determined by the Sea, Lake, and Overland Surges from Hurricanes (SLOSH) model developed by the National Oceanic and Atmospheric Administration (NOAA). In \citet{Movahednia2022a} and \citet{Movahednia2022b}, potential damages are generated by pairing flood fragility curves with spatial probability distributions of flooding derived from the Hazus tool developed by the Federal Emergency Management Agency (FEMA). In this paper, we focus on precipitation-induced fluvial flooding, a threat that is seldom studied in power system resilience literature. We use the streamflow-based forecasting method used in \citet{Austgen2022b} and \citet{Souto2022} to simulate spatially correlated damages.

%% file: modeling/_main.tex
\section{Modeling} \label{section:modeling}

\newcommand{\pathmodeling}{modeling}

\begin{figure}
    \centering
    \includegraphics[width=3in]{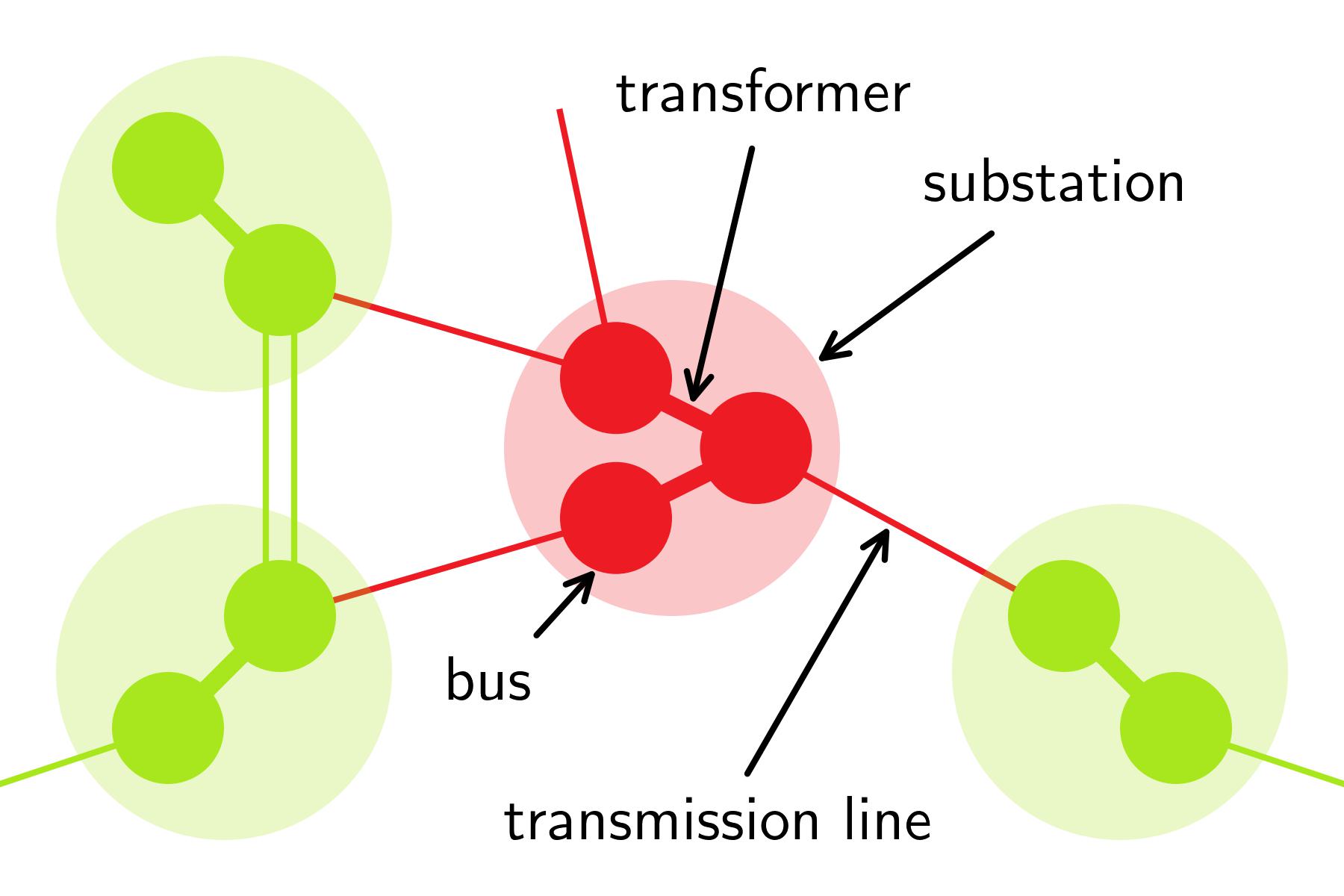}
    \caption{When a substation is flooded, its components and all adjacent transmission lines become inoperable. In this illustration, the substation shaded in light red is flooded, and the components shaded in red are consequently affected.}
    \label{fig:grid_nomenclature}
\end{figure}

For the sake of modeling power flow, we view the power grid as a graph with buses as nodes and branches (\textit{i.e.}, transmission lines and transformers) as edges. When a substation floods, we model all its buses and transformers, incident transmission lines, and associated generators as inoperable, and all its associated loads as unsatisfiable. This is illustrated in Figure~\ref{fig:grid_nomenclature}. We suppose at each substation that a discrete set of resilience levels are available for implementation by stacking temporary flood barriers like \tigerdams, and we use a power flow model to assess the load shed that results in each scenario.

\subsection{\emph{Notation}} \label{subsection:notation}
We now introduce the sets, parameters, and decision variables used in our model.

\noindent\textbf{Sets}
\begin{itemize}[align=left,leftmargin=8em,labelwidth=7em,labelsep=1em,itemsep=0em]
\item[$K$] set of substations
\item[$R = \{1, \ldots, \hat{r}\}$] set of resilience levels
\textcolor{response}{\item[$\Omega$] set of scenarios}
\item[$N, N_k$] set of buses, buses at substation $k$
\item[$E$] set of branches
\item[$N_n^+$, $N_n^-$] sets of buses neighboring bus $n$ on incoming and outgoing branches
\end{itemize}

\noindent\textbf{Parameters}
\begin{itemize}[align=left,leftmargin=8em,labelwidth=7em,labelsep=1em,itemsep=0em]
\item[$c_{kr}$] resource cost of reinforcing substation $k$ to resilience level $r$ from level $r-1$
\item[$f$] resource budget
\item[$\lambda^\text{shed}$, $\lambda^\text{over}$] unitless objective weights for load shed and overgeneration
\item[$\xi_{kr}^\omega$] flooding indicator; 1 if substation $k$ is flooded to level $r$ in scenario $\omega$; 0 otherwise
\item[$b_{nm}$] susceptance of branch $(n,m)$
\item[$\underline{p}_n^\text{gen}$, $\overline{p}_n^\text{gen}$] lower and upper bounds for power generation at bus $n$
\item[$p_n^\text{load}$] power load at bus $n$
\item[$\overline{s}_{nm}^\text{flow}$] upper bound of apparent power flow across branch $(n,m)$
\item[$n_\text{ref}$] reference bus
\item[$\overline{\theta}_\Delta$] maximum voltage phase angle difference of adjacent buses
\item[$\overline{\theta}$] maximum absolute voltage phase angle for any bus in the system
\item[$M$] arbitrarily large positive constant (for big-$M$ method)
\end{itemize}

\noindent\textbf{Decision Variables}
\begin{itemize}[align=left,leftmargin=8em,labelwidth=7em,labelsep=1em,itemsep=0em]
\item[$x_{kr} \in \{0,1\}$] mitigation indicator; 1 if substation $k$ is resilient to level $r$, 0 otherwise
\item[$\alpha_n^\omega \in \{0,1\}$] bus status indicator; 1 if bus $n$ is operational in scenario $\omega$, 0 otherwise
\item[$\beta_{nm}^\omega \in \{0,1\}$] branch status indicator; 1 if branch $(n,m)$ is operational in scenario $\omega$, 0 otherwise
\item[$\hat{p}_n^\omega, \check{p}_n^\omega \in \mathbb{R}$] power generation and overgeneration at bus $n$ in scenario $\omega$
\item[$\tilde{p}_{nm}^\omega \in \mathbb{R}$] power flow across branch $(n,m)$ in scenario $\omega$
\item[$\delta_n^\omega \in {[0,1]}$] proportion of satisfied load at bus $n$ in scenario $\omega$
\item[$\theta_n^\omega \in {[-\overline{\theta},\overline{\theta}]}$] voltage phase angle of bus $n$ in scenario $\omega$
\end{itemize}

For convenience, we denote similar types of power flow parameters and variables similarly. For example, underlining and overlining always signify lower and upper bounds, respectively (\textit{e.g.}, $\underline{p}_n^\text{gen}$ and $\overline{p}_n^\text{gen}$). Additionally, when a parameter or variable name appears in bold, it denotes the vector comprising all the indexed elements (\textit{e.g.}, $\bs{x} = [x_{kr}, \forall k \in K, \forall r \in R]$). We apply a superscript $\omega$ (\textit{e.g.}, $\bs{\xi}^\omega$) to indicate a quantity associated specifically with scenario $\omega$. All power grid parameters and variables are assumed to be in the per-unit system.

\subsection{\emph{Model}}

The discrete set of implementable resilience levels is denoted by $R$, and the decision to reinforce substation $k$ to a specific level $r$ is captured by the binary decision variable $x_{kr}$. Substation $k$ is resilient to level $r$ flooding if $x_{kr} = 1$ and is otherwise susceptible. The mitigation model comprises three sets of constraints involving these variables:
\begin{gather}
    x_{k,r+1} \le x_{k,r},\quad \forall k \in K, \forall r \in R \setminus \{\hat{r}\},
    \label{eq:con_incremental} \\
    x_{k\hat{r}} = 0,\quad \forall k \in K,
    \label{eq:con_inexorable} \\
    \sum_{k \in K} \sum_{r \in R} c_{kr} x_{kr} \le f.
    \label{eq:con_resource_hi}
\end{gather}
Constraints \eqref{eq:con_incremental} capture the cumulative nature of mitigation. We suppose the mitigation is limited either physically or practically and introduce $\hat{r}$ as the unattainable level of resilience. Flooding at or above the unattainable level is rendered inexorable by constraints \eqref{eq:con_inexorable}. Finally, we suppose the mitigation resources are limited to a budget of $f$.
Supposing a marginal cost $c_{kr}$ associated with each decision $x_{kr}$, this is captured by the binary knapsack constraint \eqref{eq:con_resource_hi}. For brevity, we hereafter refer to the constraints of the mitigation decision making problem as
\begin{equation}
    \mathcal{X} = \{\boldsymbol{x} \in \{0,1\}^{|K \times R|} :
                    \eqref{eq:con_incremental},
                    \eqref{eq:con_inexorable},
                    \eqref{eq:con_resource_hi}\}. \label{eq:resilience_X}
\end{equation}

The DC power flow model that we use in the recourse problem is based on three simplifying assumptions: (1) branch conductance is negligible relative to susceptance and may be ignored, (2) bus voltage magnitudes are approximately one per unit, and (3) the difference between voltage phase angles of adjacent buses is small such that the sine of that difference is approximately linear and the cosine is approximately one \citep{Molzahn2019}.
These assumptions lead to the reactive power flows being zero and the active power flows obeying a linear relationship with the bus voltage phase angles.
The model may be viewed as an extension of a capacitated network flow problem with multiple sources and sinks. Importantly, the solution space is additionally confined by the complicating Ohm's Law constraints.
Our adaptation of the DC power flow model is as follows:
\begin{mini!}[2]<b>
    {}{\sum_{n \in N} \lambda^\text{shed} p_n^\text{load} (1 - \delta_n^\omega) + \lambda^\text{over} \check{p}_n^\omega}{\label{eq:loss}}
    {\mathcal{L}(\boldsymbol{x},\boldsymbol{\xi}^\omega) =} \label{eq:dc_min_load_shed}
    \addConstraint{\alpha_n^\omega = \prod_{r \in R} \left(1 - \xi_{kr}^\omega \left(1 - x_{kr}\right)\right),}
                  {\quad}{\forall k \in K, \forall n \in N_k,} \label{eq:dc_def_alpha}
    \addConstraint{\beta_{nm}^\omega = \alpha_n^\omega \alpha_m^\omega,}
                  {\quad}{\forall (n,m) \in E,} \label{eq:dc_def_beta}
    \addConstraint{\hat{p}_n^\omega - \check{p}_n^\omega - p_n^\text{load} \delta_n^\omega
                   + \sum_{m \in N_n^+} \tilde{p}_{mn}^\omega
                   - \sum_{m \in N_n^-} \tilde{p}_{nm}^\omega
                   = 0,}
                   {\quad}{\forall n \in N,} \label{eq:dc_kcl}              
    \addConstraint{M \left(\beta_{nm}^\omega - 1\right) \le -\tilde{p}_{nm}^\omega - b_{nm} (\theta_n^\omega - \theta_m^\omega),}
                  {\quad}{\forall (n,m) \in E,} \label{eq:dc_ohms_law_le}
    \addConstraint{M \left(1 - \beta_{nm}^\omega\right) \ge -\tilde{p}_{nm}^\omega - b_{nm} (\theta_n^\omega - \theta_m^\omega),}
                  {\quad}{\forall (n,m) \in E,} \label{eq:dc_ohms_law_ge}
    \addConstraint{\theta_n^\omega - \theta_m^\omega \ge -2 (1 - \beta_{nm}^\omega) \overline{\theta} - \beta_{nm}^\omega \overline{\theta}_\Delta,}
                  {\quad}{\forall (n, m) \in E,} \label{eq:phase_angle_diff_le}
    \addConstraint{\theta_n^\omega - \theta_m^\omega \le 2 (1 - \beta_{nm}^\omega) \overline{\theta} + \beta_{nm}^\omega \overline{\theta}_\Delta,}
                  {\quad}{\forall (n, m) \in E,} \label{eq:phase_angle_diff_ge}
    \addConstraint{-\overline{s}_{nm}^\text{flow} \beta_{nm}^\omega \le \tilde{p}_{nm}^\omega \le \overline{s}_{nm}^\text{flow} \beta_{nm}^\omega,}
                  {}{\forall (n, m) \in E,} \label{eq:dc_flow_limits}
    \addConstraint{\underline{p}_n^\text{gen} \alpha_n^\omega \le \hat{p}_n^\omega \le \overline{p}_n^\text{gen} \alpha_n^\omega,}
                  {\quad}{\forall n \in N,} \label{eq:dc_generation_limits}
    \addConstraint{0 \le \check{p}_n^\omega \le \hat{p}_n^\omega,}
                  {\quad}{\forall n \in N,} \label{eq:dc_overgeneration_limit}
    \addConstraint{\theta_{n_\text{ref}}^\omega = 0}
                  {}{} \label{eq:dc_ref_voltage_phase_angle}.
\end{mini!}
Here, variables are also constrained as specified in Section~\ref{subsection:notation}. The objective \eqref{eq:dc_min_load_shed} is to minimize the weighted combination of load shed and overgeneration. Constraints \eqref{eq:dc_def_alpha} and \eqref{eq:dc_def_beta} relate the operational statuses of buses and branches to those of the substations. Note that these are equality constraints. As such, every component's operational status is perfectly determined by first-stage decisions $\boldsymbol{x}$ and flooding uncertainty realization $\boldsymbol{\xi}^\omega$.
Though these constraints are nonlinear in their presented logical form, they admit linear reformulations:
\begin{alignat}{2}
    & \alpha_n^\omega \ge \displaystyle \sum_{r \in R} \left(1 - \xi_{kr}^\omega \left(1 - x_{kr}\right)\right) - |R| + 1,
    && \quad \forall k \in K, \forall n \in N_k, \\
    & \alpha_n^\omega \le 1 - \xi_{kr}^\omega \left(1 - x_{kr}\right),
    && \quad \forall k \in K, \forall n \in N_k, \forall r \in R, \\
    & \beta_{nm}^\omega \ge \alpha_n^\omega + \alpha_m^\omega - 1, && \quad \forall (n,m) \in E, \\
    & \beta_{nm}^\omega \le \alpha_n^\omega,                && \quad \forall (n,m) \in E, \\
    & \beta_{nm}^\omega \le \alpha_m^\omega,                && \quad \forall (n,m) \in E.
\end{alignat}
Kirchhoff's Current Law (KCL), the power flow equivalent of flow balance, is imposed by constraints \eqref{eq:dc_kcl}. In the standard DC power flow formulation, Ohm's Law is represented as the equality constraint ${p_{nm}^\omega = -b_{nm} (\theta_n^\omega - \theta_m^\omega)}$. To ensure out-of-service branches are treated as open circuits, we employ the big-$M$ technique used in \citet{Coffrin2011} in constraints \eqref{eq:dc_ohms_law_le} and \eqref{eq:dc_ohms_law_ge} to enforce the equality only for operational branches. Similarly, constraints \eqref{eq:phase_angle_diff_le} and \eqref{eq:phase_angle_diff_ge} impose limits on the differences of phase angles only for buses joined by one or more operational branches. Constraints \eqref{eq:dc_flow_limits} impose conditional lower and upper bounds on power flows. Power generation lower and upper bounds are imposed by constraints \eqref{eq:dc_generation_limits} and are dependent on the operational status of the corresponding bus. Constraints \eqref{eq:dc_overgeneration_limit} ensure no more power is overgenerated than is generated at each generator, and constraint \eqref{eq:dc_ref_voltage_phase_angle} ensures the voltage phase angle of the reference bus is exactly zero.

The role of overgeneration in the recourse model is to allow generators to effectively operate below their lower limit. While the objective's primary purpose is to minimize load shed, its serves the secondary purpose of minimizing net violations of this soft constraint. Modeling generation this way admits relatively complete recourse -- the solution with $\hat{p}_n^\omega = \check{p}_n^\omega = \underline{p}_n^\text{gen}$ for all $n \in N$ and all other power flow variables equal to zero is always feasible.

\begin{figure}
    \centering
    \includegraphics[width=6.2in]{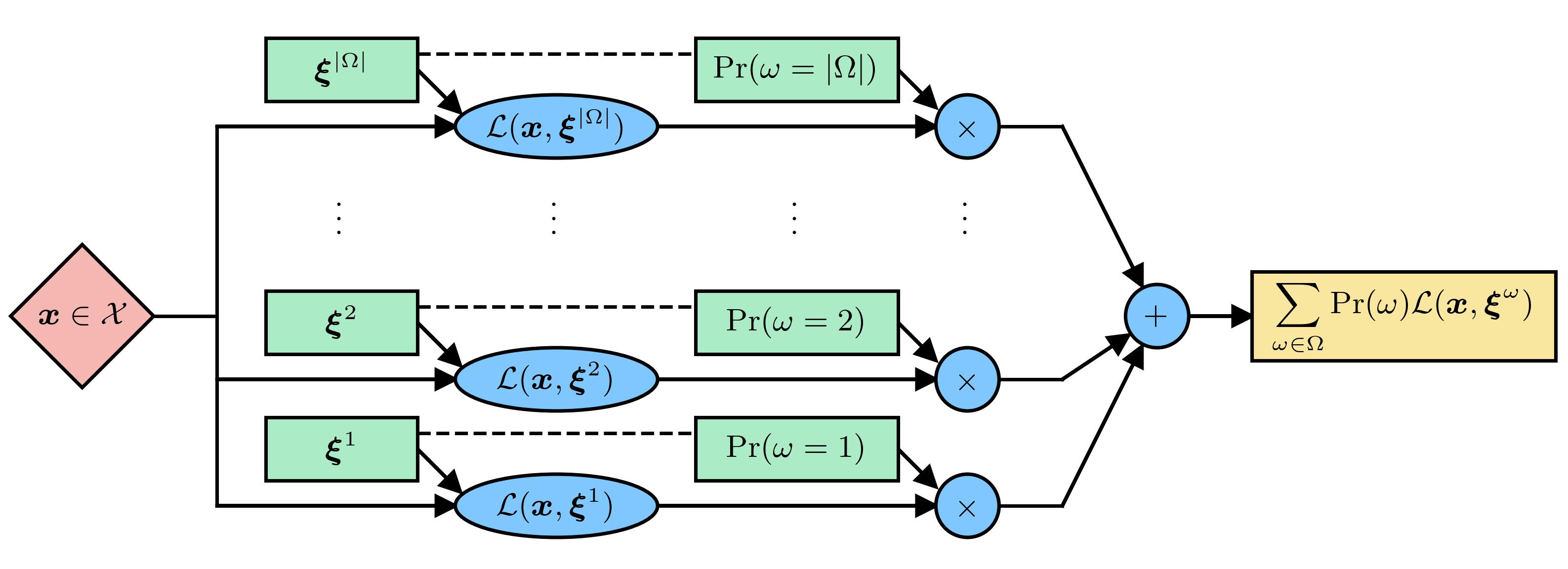}
    \caption{\textcolor{response}{A diagram of our two-stage model. First, a decision maker commits to a specific mitigation solution $\bs{x} \in \mathcal{X}$. For each scenario $\omega \in \Omega$, the load shed-minimizing grid operational response is computed via recourse problem $\mathcal{L}(\bs{x}, \bs{\xi}^\omega)$, a physics-based model that captures power flow limitations following mitigation solution $\bs{x}$ and flooding realization $\bs{\xi}^\omega$. The efficacy of the mitigation solution is measured by the probability-weighted sum of loss (\textit{i.e.}, the expected weighted combination of load shed and overgeneration) across all sample scenarios. The goal of our model is to identify the mitigation solution (red diamond) that minimizes the resulting expected load shed (yellow rectangle).}}
    \label{fig:model_diagram}
\end{figure}

With $\mathcal{X}$ and $\mathcal{L}$ as defined in \eqref{eq:resilience_X} and \eqref{eq:loss}, respectively, we formulate our two-stage stochastic programming (SP) model as
\begin{equation}
    \min_{\bs{x} \in \mathcal{X}} \sum_{\omega \in \Omega} \prob(\omega) \mathcal{L}(\bs{x}, \bs{\xi}^\omega)
    \tag{SP} \label{eq:SP}.
\end{equation}
In this model, mitigation decisions are evaluated by how well the power grid is able to perform in the aftermath of the hurricane given those decisions. Based on our assumption that a hurricane is naturally (\textit{i.e.}, randomly) occurring, we model each scenario $\omega$ in the ensemble $\Omega$ is having probability $\prob(\omega)$ of occurrence and weight the associated loss accordingly in the objective function.
\textcolor{response}{We illustrate how the expected loss resulting from a specific mitigation solution is computed in our two-stage stochastic programming formulation in Figure~\ref{fig:model_diagram}.}

Given the prevalence of this unfettered ``minimize load shed'' approach in power system resilience research \citep{Coffrin2011,Pierre2018,Moreno2020,Movahednia2022a}, we do the same. This objective tries to minimize the system-wide load shed but may leave certain loads completely unsatisfied in some or all of the scenarios. If the decision maker has goals of satisfying individual loads as much as possible, an alternative solution could be to introduce a service level constraint for each load bus, \textit{i.e.},
\begin{equation}
    \sum_{\omega \in \Omega} \textrm{Pr}(\omega) \delta_n^\omega \ge \underline{\delta}_n, \quad \forall n \in N.
    \label{eq:service_level_naive}
\end{equation}
where $\underline{\delta}_n$ is the expected load satisfaction level imposed at bus $n$. Such constraints are often used to inform the management of variable renewable electricity generation \citep{Bienstock2014,Roald2018} and less commonly used in the resilience literature \citep{Luo2016}. With or without such constraints, in this application, load shed occurs mainly because bus outages prevent loads from being served and branch outages diminish the capabilities of the system to flow power to where it is demanded at in-service buses. For given $\boldsymbol{x}$ and $\boldsymbol{\xi}^\omega$, most loads are either completely satisfied or completely unsatisfied. Due to the nearly discrete nature of load satisfaction, and the discrete nature of scenarios' impact on substations, the optimal solutions will show unintuitive and non-monotonic changes as the budget is increased as is typical of discrete knapsack-type problems. We highlight this phenomenon in Section~\ref{subsection:solution_analysis} while presenting results from experiments without the service level constraints. Our experiments with service level constraints of type \eqref{eq:service_level_naive}, though not presented here, show similar behavior.

%% file: casestudies/_main.tex
\section{Case Study Development} \label{section:case_studies}

In this section, we describe our development of the two case studies to which we apply our model.
The cases are both built on the same power grid instance.
However, they are each based on different flooding events and additionally differ in the methodology used to model the flooding.
We take special care to ensure a high degree of geographic realism in the case studies.

\subsection{\emph{Power Grid}}

Both of our case studies are based on the ACTIVS 2000-bus synthetic grid \citep{Birchfield2017, Gegner2016}. Though the grid is synthetic, it is designed to be statistically similar to the Texas Interconnection, and its 1250 substations are geographically defined in that region. For our application, using the original data presents two problems. First, some of the synthesized coordinates for the substations either reside in or are in close proximity to major water bodies such that they are submerged under normal circumstances, at least according to the tools we use to model flooding. Second, though the grid fits the geographical scale needed for our application, embedding multiple 2000-bus power flow instances in a stochastic program is computationally limiting.

\subsubsection{\emph{Coordinate Remapping}}

To augment the geographical realism in our case study, we relocate all 1250 substations from the original ACTIVS 2000-bus grid to locations of actual substations in the state of Texas. We achieve this by computing a minimum-distance mapping of the 1250 original substation coordinates to the coordinates of substations from the Homeland Infrastructure Foundation-Level Data (HIFLD) Electric Substations dataset \citep{HIFLD_Electric_Substations} which contains information about real-world substations across the U.S. states and territories.

We let $A$ denote the set of substations from the ACTIVS 2000-bus dataset and $B$ the subset of substations from the HIFLD Electric Substations dataset that are located in Texas. Using $c_{ab}$ to denote the distance from $a \in A$ to $b \in B$ and $x_{ab}$ the decision to map $a$ to $b$, we compute the mapping by solving the unbalanced assignment problem
\begin{mini!}
{}{\sum_{a \in A} \sum_{b \in B} c_{ab} x_{ab}}{}{} \label{eq:assignment_min_distance}
\addConstraint{\sum_{b \in B} x_{ab} \ge 1,}{\quad}{\forall a \in A} \label{eq:assignment_every_a}
\addConstraint{\sum_{a \in A} x_{ab} \le 1,}{\quad}{\forall b \in B} \label{eq:assignment_at_most_one_b}
\addConstraint{x_{ab} \in \{0, 1\},}{\quad}{\forall a \in A, \forall b \in B} \label{eq:assignment_integrality}.
\end{mini!}

Objective \eqref{eq:assignment_min_distance} minimizes the cost (\textit{i.e.}, total distance) of the mapping. Constraints \eqref{eq:assignment_every_a} ensure each $a \in A$ is mapped to some $b \in B$. Constraints \eqref{eq:assignment_at_most_one_b} ensure no more than one $a \in A$ is mapped to any $b \in B$. Finally, \eqref{eq:assignment_integrality} simply imposes integrality. However, because the constraint matrix formed by \eqref{eq:assignment_every_a} and \eqref{eq:assignment_at_most_one_b} is totally unimodular, the integrality constraints may be relaxed so that each $x_{ab}$ resides in the unit interval.

\subsubsection{\emph{Reduction}}

\begin{figure}
    \centering
    \setlength{\belowcaptionskip}{\baselineskip}
    \begin{subfigure}{0.49\textwidth}
        \includegraphics[width=\textwidth]{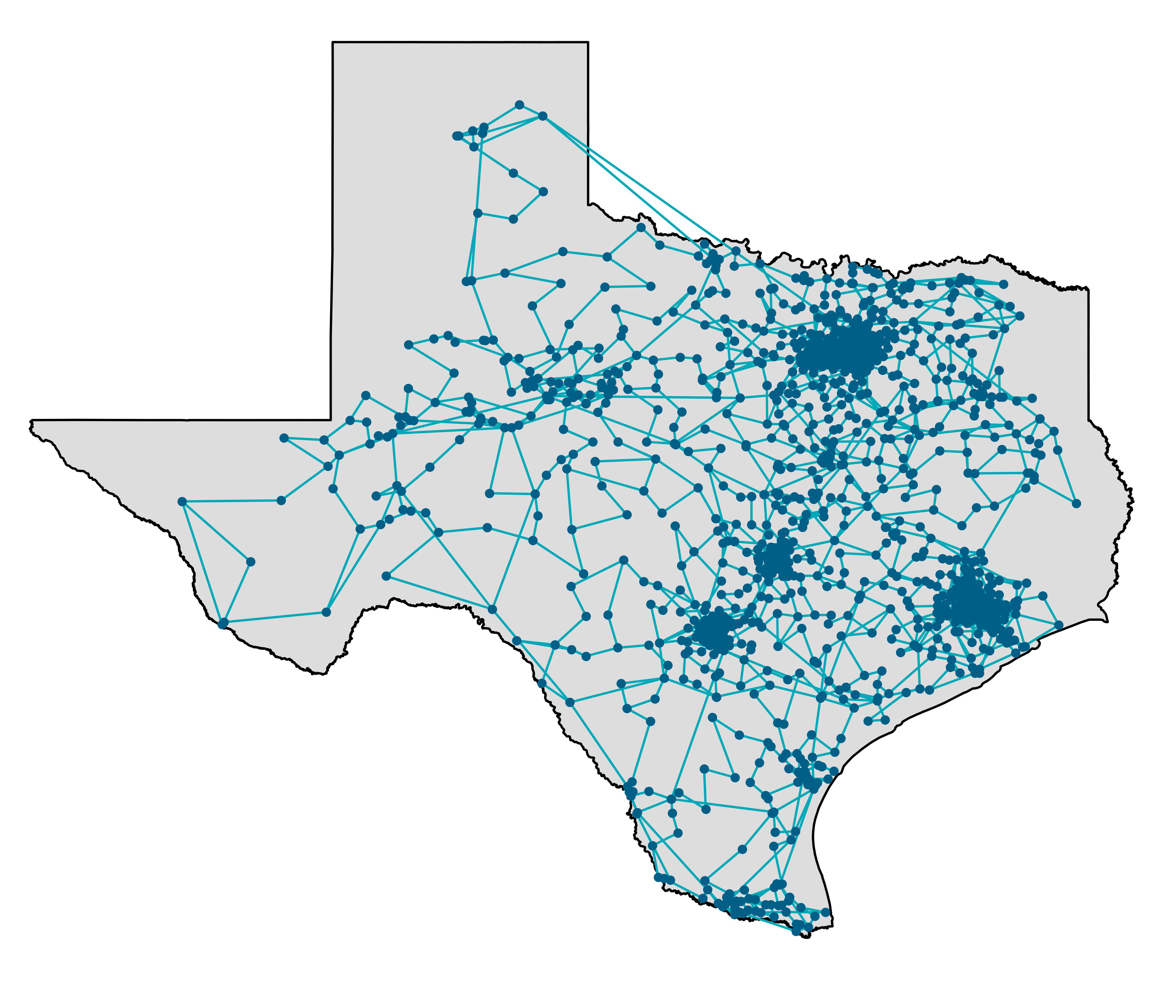}
        \caption{}
        \label{fig:grid_before_reduction}
    \end{subfigure}
    \begin{subfigure}{0.49\textwidth}
        \includegraphics[width=\textwidth]{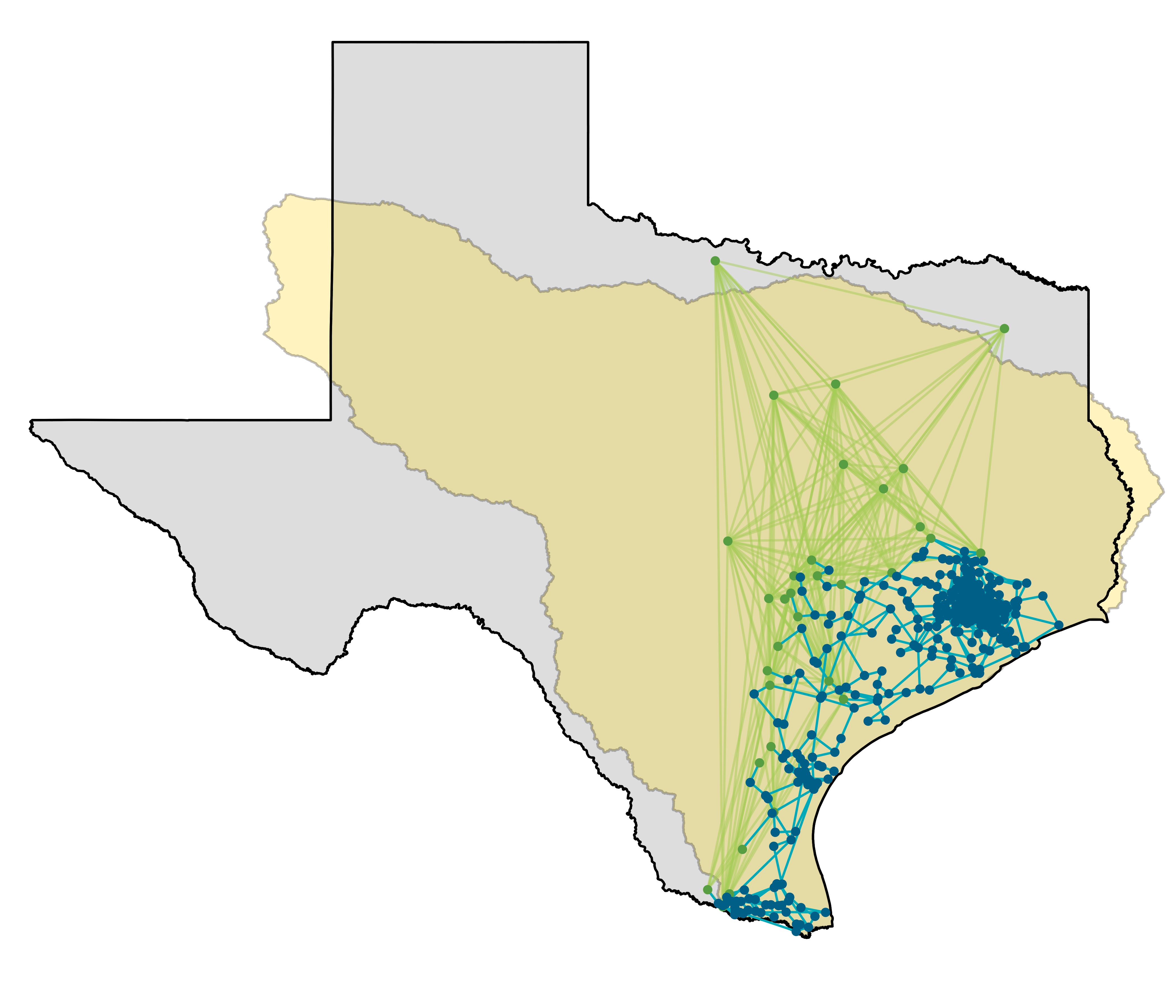}
        \caption{}
        \label{fig:grid_after_reduction}
    \end{subfigure}
    \caption{Comparing (a) the full original grid with (b) the reduced grid, we observe the number of buses and branches are reduced by more than 50\%. In (b), teal represents the parts retained from the original grid, and green represents the parts introduced by the reduction. The Texas-Gulf region, the area modeled in our streamflow simulations, is shaded yellow.}
    \label{fig:texas_maps}
\end{figure}

To reduce the size of the optimization problem, we compute a network reduction of the ACTIVS 2000-bus synthetic grid using the electrical equivalent (EEQV) feature from {PSS\textregistered E} \citep{PSSE2013}. The reduced network aggregates buses in the inland region that are not exposed to flooding conditions but retains the parts of the grid in the coastal region that are affected by flooding. Table~\ref{tab:grid_before_and_after} details the high-level effects of the reduction and Figures~\ref{fig:grid_before_reduction}~and~\ref{fig:grid_after_reduction} illustrate the geographical differences.

\begin{table}
    \centering
    \caption{Power grid characteristics before and after the reduction was performed.}
    \begin{tabular}{|l|r|r|}
        \hline
        Grid Characteristic & Before & After \\
        \hline
        Substations (\#) & 1250 & 362 \\
        Buses (\#) & 2000 & 663 \\
        Transformers (\#) & 860 & 369 \\
        Transmission Lines (\#) & 2346 & 1140 \\
        Generators (\#) & 544 & 254 \\
        \hline
        Generation Capacity (GW) & 96.2915 & 50.9779 \\
        Load (GW) & 67.1092 & 39.6860 \\
        \hline
    \end{tabular}
    \label{tab:grid_before_and_after}
\end{table}

\subsection{\emph{Mitigation}}
\begin{figure}
    \centering
    \ifdefined\epsswitch
        \includegraphics[width=3in]{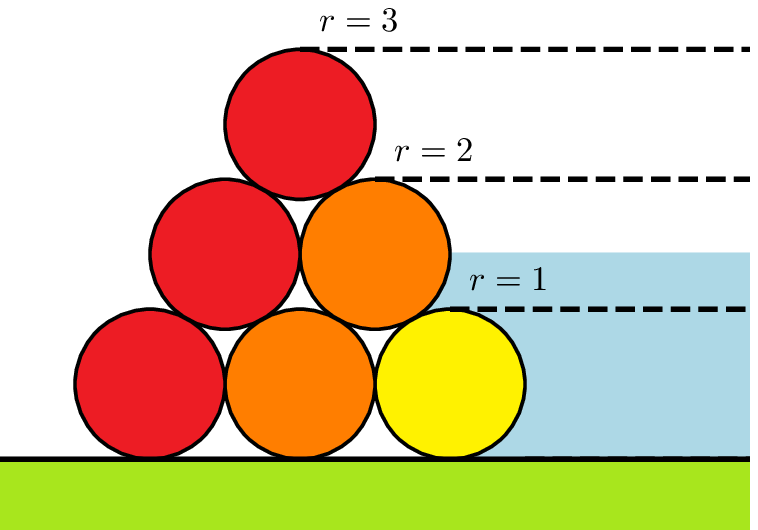}
    \else
        \includegraphics[width=3in]{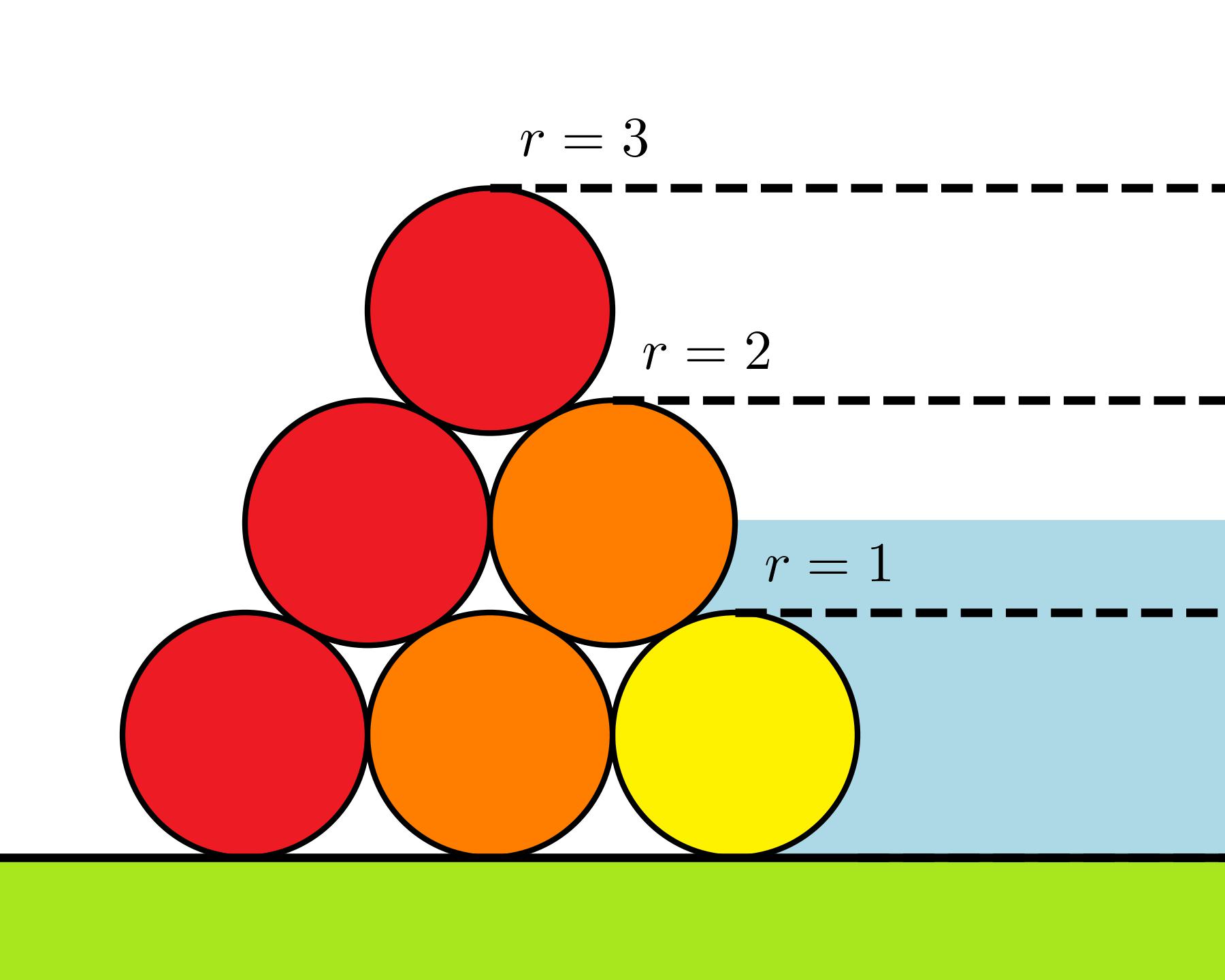}
    \fi
    \caption{The cross-sectional view of a stacked \tigerdam~flood barrier. The colors indicate the resources required to attain the next-highest level of resilience. Though level $r=3$ is illustrated here, we model this level as being unattainable (\textit{i.e.}, $\hat{r} = 3$) in our case studies.}
    \label{fig:mitigation_levels}
\end{figure}

In our case studies, we let the set of implementable resilience levels be $R = \{1, 2, 3\}$ and the unattainable level of resilience be $\hat{r} = 3$. This $R$ represents the discrete set of resilience levels that may be achieved by extending and stacking \tigerdams~around a substation as shown in Figure~\ref{fig:mitigation_levels}. The figure illustrates the cross section of a modular flood barrier. Level $r = 0$ represents the case in which no mitigation is implemented such that any level of flooding renders the substation inoperable.

The assumed cross-sectional radius of a \tigerdam~is 0.268 meters such that level $r = 1$ resilience protects against 0.534 meters of flooding and level $r = 2$ protects against exactly 1 meter of flooding. Level $r = 3$ resilience would protect against 1.464 meters of flooding if its implementation were to be allowed.

The number of \tigerdams~required to make a substation level $r$ resilient depends on both the level $r$ and the substation's perimeter. The ACTIVS 2000-bus grid has components defined on four different voltage levels: 115~kV, 161~kV, 230~kV, and 500~kV. We assume a single \tigerdam~segment is sufficient to cover the perimeter of substations for which the highest-voltage component is either 115~kV or 161~kV. If the highest-voltage component is rather 230~kV or 500~kV, then either 2 or 3 segments, respectively, are assumed to be needed. Figure~\ref{fig:mitigation_levels} illustrates the marginal costs of resilience. Because stacking requires additional \tigerdams~to be placed in the base layers, the marginal cost is linearly increasing in $r$. The number of marginal resources required to implement each level $r \in R$ for each size of substation is summarized in Table~\ref{tab:mitigation_costs}.

\begin{table}
    \centering
    \caption{Marginal resource requirements for protecting differently sized substations.}
    \begin{tabular}{|r|ccc|}
        \hline
        Highest-Voltage Component & \multicolumn{3}{c|}{Resilience Level} \\
        \hline
        & $r=1$ & $r=2$ & $r=3$ \\
        \hline
        115 kV or 161 kV & 1 & 2 & 3 \\
                  230 kV & 2 & 4 & 6 \\
                  500 kV & 3 & 6 & 9 \\
        \hline
    \end{tabular}
    \label{tab:mitigation_costs}
\end{table}

In subsequent sections, we use the terms ``resilience level'' and ``flood level'' frequently. The terms are not interchangeable, but they are related. The term ``flood level $r$'' indicates flooding that requires implementing at least resilience level $r$  to mitigate. Take Figure~\ref{fig:mitigation_levels} for example. Level $r = 1$ resilience is insufficient for mitigating the flooding shown in light blue. However, level $r = 2$ resilience is sufficient and not excessive. Thus, the illustrated flooding is level $r = 2$.

\subsection{\emph{Flooding Scenarios}}
To produce hurricane-induced flooding scenarios, we apply the National Water Model (NWM) to data available for historical events prior to their landfall. The two historical events we study are Tropical Storm Imelda from 2019 and Hurricane Harvey from 2017, respectively the 68th and 2nd most costly U.S. natural disasters since 1980 \citep{BillionDollarDisasters2022}. These events both affected the Texas coastal region, especially Houston, and were noteworthy for the rainfall they produced \citep{Blake2018,Latto2020}.

The National Oceanographic and Atmospheric Administration (NOAA) launched the NWM in 2016 as part of an effort to improve water analysis and prediction capabilities and promote resilience to water risk \citep{NWM2016}. Broadly, the model operates on meteorological inputs like temperature and precipitation to produce hydrological outputs such as stream flow. More specifically, the NWM is a set of configurations of the community-based Weather Research and Forecasting Hydrologic (WRF-Hydro) framework developed by the National Center for Atmospheric Research (NCAR) \citep{Gochis2020}. For the continental United States (CONUS), NWM uses three configurations for short-, medium-, and long-range forecasting, but only short- and medium-range forecasting are practical for our application. The Short-Range Forecast (SRF) is a single deterministic forecast out to 18 hours and the Medium-Range Forecast (MRF) is a 7-member ensemble forecast with the longest ensemble member capturing out to 10 days. A principal difference is that SRF is driven by forcing from the High-Resolution Rapid Refresh (HRRR) model whereas MRF is driven by the Global Forecasting System (GFS).

Because we target storms that primarily affected Texas, we opt not to model all of CONUS in our forecasting. Rather, we target only the Texas-Gulf region. In the Hydrological Unit Code (HUC) taxonomy, the Texas-Gulf region corresponds to two-digit HUC \#12 \citep{USGS2022a}. This region is mostly contained in the state of Texas but also includes small sections of Louisiana and New Mexico as shown in Figure~\ref{fig:texas_maps}.

Importantly, the methodologies we develop for forecasting Tropical Storm Imelda and Hurricane Harvey account only for fluvial flooding or flooding caused by streams and rivers overflowing their banks. They do not account for pluvial flooding, also known as ponding, or flooding caused by tidal surge. Even so, we do not model the uncertainty of fluvial flooding the same way for both storms. Rather, we rely on different data and develop different methodologies for the two. We present more details about our approaches hereafter.

\subsubsection{\emph{Tropical Storm Imelda}}

To produce forecasts for Tropical Storm Imelda, we follow the approach of \citet{Wu2022}. This approach accounts for the uncertainty of the hurricane's path and precipitation intensity, and it incorporates ideas from both the SRF and MRF products. The SRF product provides better spatiotemporal resolution and thus presumably better forecasts than the MRF product, but its simulated time interval of 18 hours is too short for our application. Decision makers would want as much time as possible prior to landfall to implement any short-term resilience measures; however, if the SRF is produced too soon in advance of actual landfall then the simulated interval would not even include the landfall event much less the full extent of the flooding inundation. The 10-day outlook provided by the MRF is more suitable.

To remedy these issues, our approach perturbs the GFS weather inputs to the MRF according to the HRRR weather inputs that ordinarily drive the SRF. Tropical Storm Imelda was a particularly short-formed and short-lived hurricane \citep{Blake2018}. It formed on September 17, 2019 at 12:00 UTC, and the National Hurricane Center (NHC) issued the first advisory 5 hours later \citep{NHC2019a}. For the GFS forecast initialized at 12:00 UTC and the HRRR forecasts initialized at 12:00 UTC, 13:00 UTC, 14:00 UTC, and 15:00 UTC, we assess the spatial maximum of accumulated precipitation over their shared 15-hour interval between September 17, 2019 at 15:00 UTC and September 18, 2019 at 06:00 UTC. We also assess the coordinates of the storm's landfall location which is assumed to be its projected location on September 17, 2019 at 22:00 UTC, the time it made landfall in the GFS forecast. These assessments are summarized in Table~\ref{tab:weather_forecast_summary}.

\begin{table}
    \centering
    \caption{Properties of the GFS and HRRR forecasts of Imelda used to parameterize the uncertainty.}
    \begin{tabular}{|c|c|c|c|}
        \hline
         Weather       & Initialization & Spatial Max. of 15-hour    & Landfall \\
         Forecast Tool & Time           & Accumulated Precipitation  & Coordinates \\
         \hline
         GFS  & 12:00 UTC & 1.6865 mm & 95.28\textdegree W, 29.02 \textdegree N \\
         HRRR & 12:00 UTC & 5.1004 mm & 95.17\textdegree W, 29.25 \textdegree N \\
         HRRR & 13:00 UTC & 3.5985 mm & 95.21\textdegree W, 29.10 \textdegree N \\
         HRRR & 14:00 UTC & 2.7876 mm & 95.31\textdegree W, 29.08 \textdegree N \\
         HRRR & 15:00 UTC & 3.3803 mm & 95.22\textdegree W, 29.09 \textdegree N \\
         \hline
    \end{tabular}
    \label{tab:weather_forecast_summary}
\end{table}

\begin{figure}
    \centering
    \includegraphics{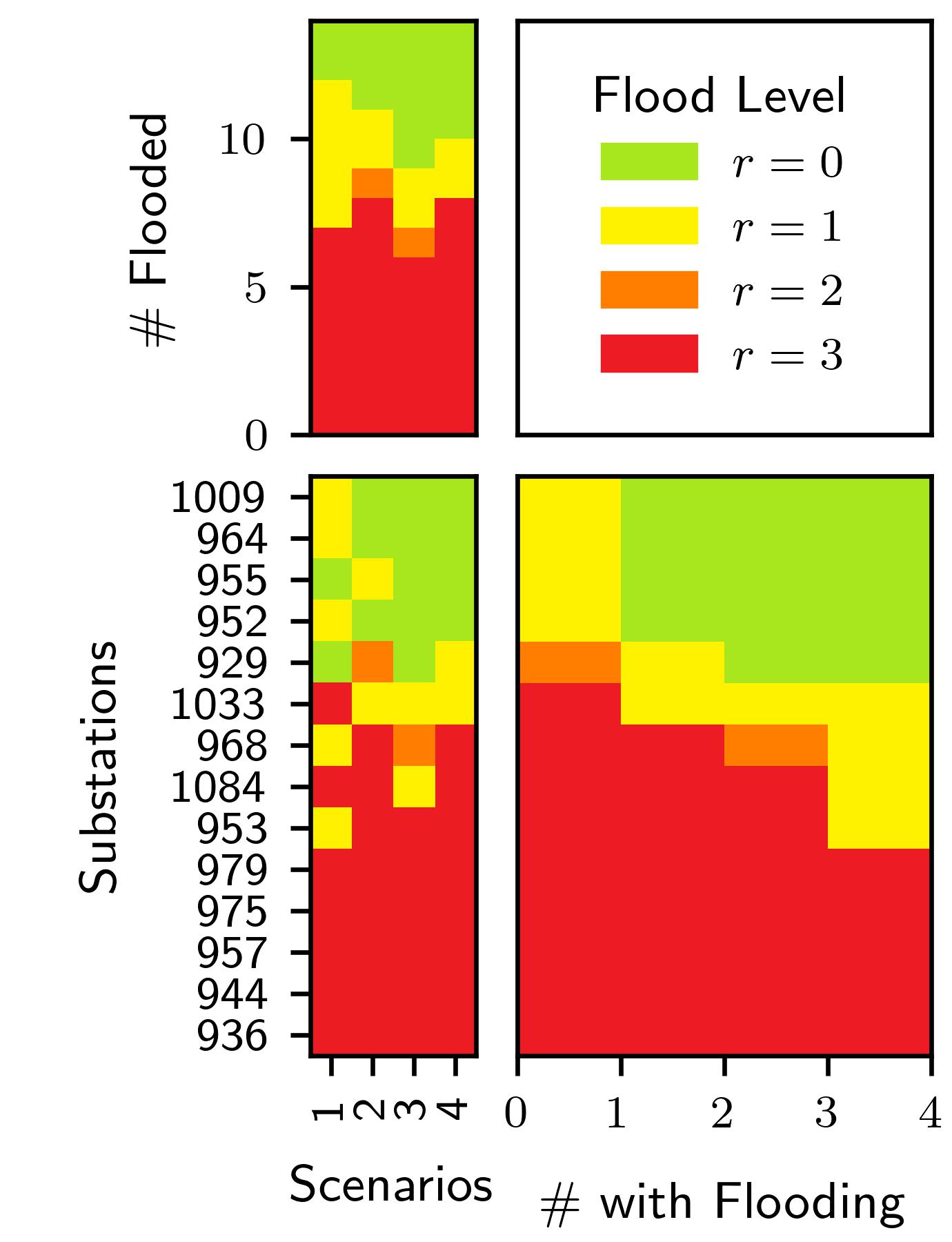}
    \caption{Tropical Storm Imelda scenario flood levels by scenario and by substation.}
    \label{fig:imelda_uncertainty_heatmap}
\end{figure}

In every case, the HRRR forecasts place the storm's landfall location northeast of its location in the GFS forecast. The HRRR forecasts also all yield a greater spatial maximum of accumulated precipitation. To produce flooding scenarios, we manipulate the GFS forecast to match each HRRR forecast in two ways: first by scaling the precipitation to match the spatial maximum of accumulated precipitation and second by spatially translating all parameters according to the difference in landfall coordinates. These forecasts are then fed to WRF-Hydro configured to run as the MRF in every way except for the input perturbations. For each scenario, we take the flood level at each substation to be the temporal maximum of the flood levels from the forecast. An aggregate spatial perspective of the resulting flooding scenarios is presented in Figure~\ref{fig:imelda_uncertainty_map} and an illustration of the \textcolor{response}{flooding in each scenario} in Figure~\ref{fig:imelda_uncertainty_heatmap}.
\textcolor{response}{Though this approach to generating scenarios is not probabilistic, we suppose the forecasts are equiprobable in the parameterization of our two-stage decision model.}

\subsubsection{\emph{Hurricane Harvey}}

Our approach to forecasting Tropical Storm Imelda is designed so that the insights gained from the high-resolution weather forecast may be applied to the medium-resolution weather forecast used to drive the more suitable medium-range forecasting tool. For those insights to be meaningful, the forecasts must be initialized late enough to simulate through the time of landfall; however, this restriction leaves decision makers with only about 15 hours to form and implement a mitigation strategy. \textcolor{response}{Moreover, our approach to forecasting Imelda does not provide a justifiable basis for assigning probabilities to scenarios. For forecasting Hurricane Harvey, we improved on these shortcomings by developing an approach based on that from \citet{Kim2021} that incorporates probability and permits more time for mitigation planning and execution.}

The approach considers a single uncertain parameter, the storm's path, and it does so by accounting for the storm's cone of uncertainty. The cone of uncertainty is the region around the projected storm path that is believed to contain the center of the storm with two-thirds probability based on NOAA's forecasting errors in the basin over the past five years \citep{NHC2022}. Forecast errors increase in the outlook time, hence why it is a \textit{cone} of uncertainty -- the region expands as the outlook time increases.

\begin{figure}
    \centering
    \ifdefined\epsswitch
        \includegraphics[width=3in]{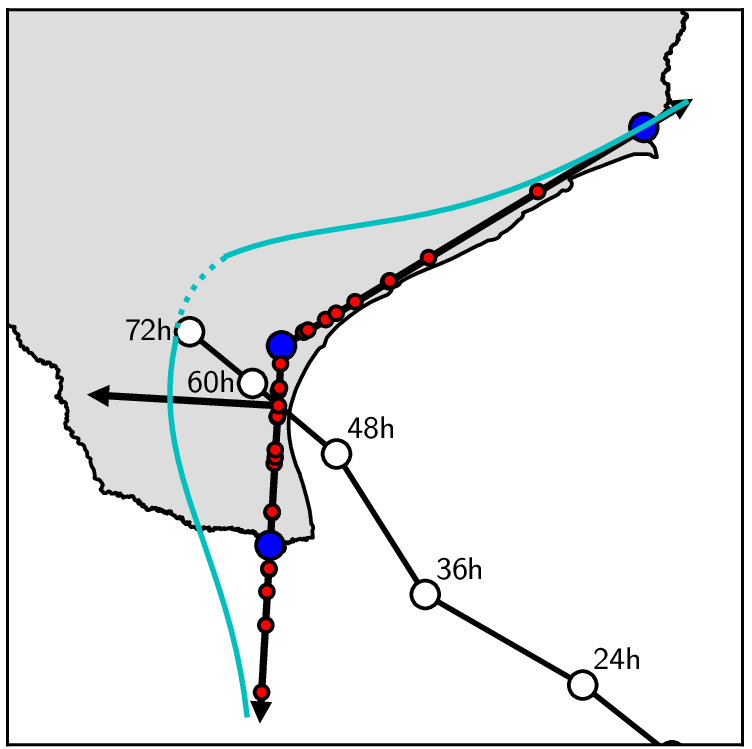}
    \else
        \includegraphics[width=3in]{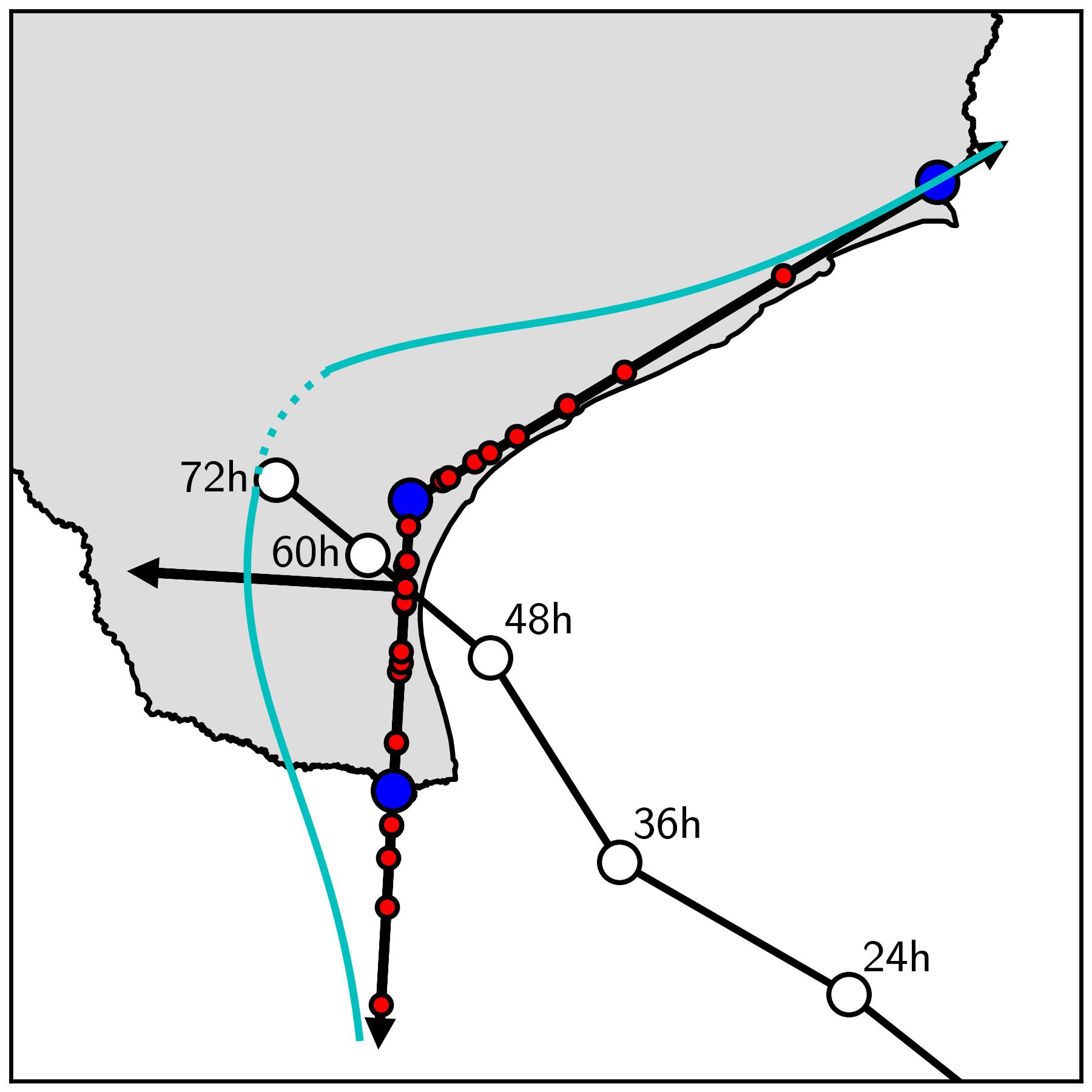}
    \fi
    \caption{The distribution of Hurricane Harvey's landfall is defined on a piece-wise linear approximation of the coastline. The blue markers define the coastline approximation, the white markers are the path from the GFS forecast, and the red markers are sample landfall locations.}
    \label{fig:harvey_landfall_distribution}
\end{figure}

A goal of this approach is to produce plausible scenarios around 48 hours ahead of the hurricane's projected landfall. Because NOAA's advisories are not necessarily timed in this way, we leverage an advisory that projected Harvey to make landfall roughly 57 hours from the time the advisory was issued, and the radius of the cone of uncertainty with that outlook time in the Atlantic basin in 2017 is approximately 89 nautical miles. Using this information, we project a normal distribution of landfall location onto a piece-wise linear approximation of the Texas coastline. The distribution is parameterized such that the storm makes landfall within 89 nautical miles of the projected (\textit{i.e.}, mean) landfall location with two-thirds probability. This distribution is illustrated in Figure~\ref{fig:harvey_landfall_distribution}.

With the distribution constructed, we sample 25 landfall locations using a stratified approach. To produce a flooding scenario, we spatially translate the GFS data from August 24, 2017, at 00:00 UTC so that its landfall location matches that of the corresponding sample, and then we feed the data to WRF-Hydro configured as MRF. This is similar to our approach for Imelda except that we do not also scale the precipitation data.

\begin{figure}
    \centering
    \begin{subfigure}{0.49\textwidth}
        \centering
        \includegraphics[width=3in]{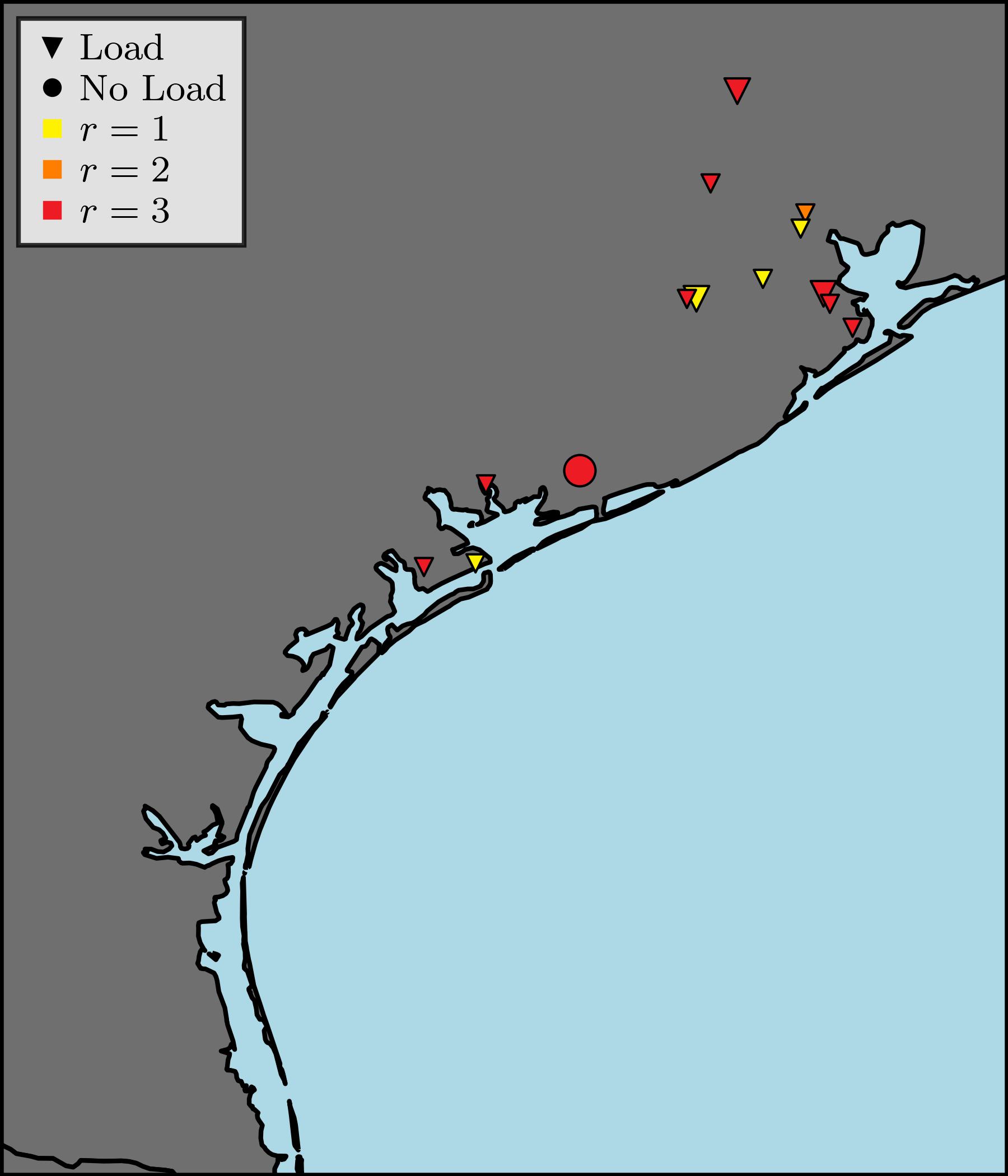}
        \caption{}
        \label{fig:imelda_uncertainty_map}
    \end{subfigure}
    \begin{subfigure}{0.49\textwidth}
        \includegraphics[width=3in]{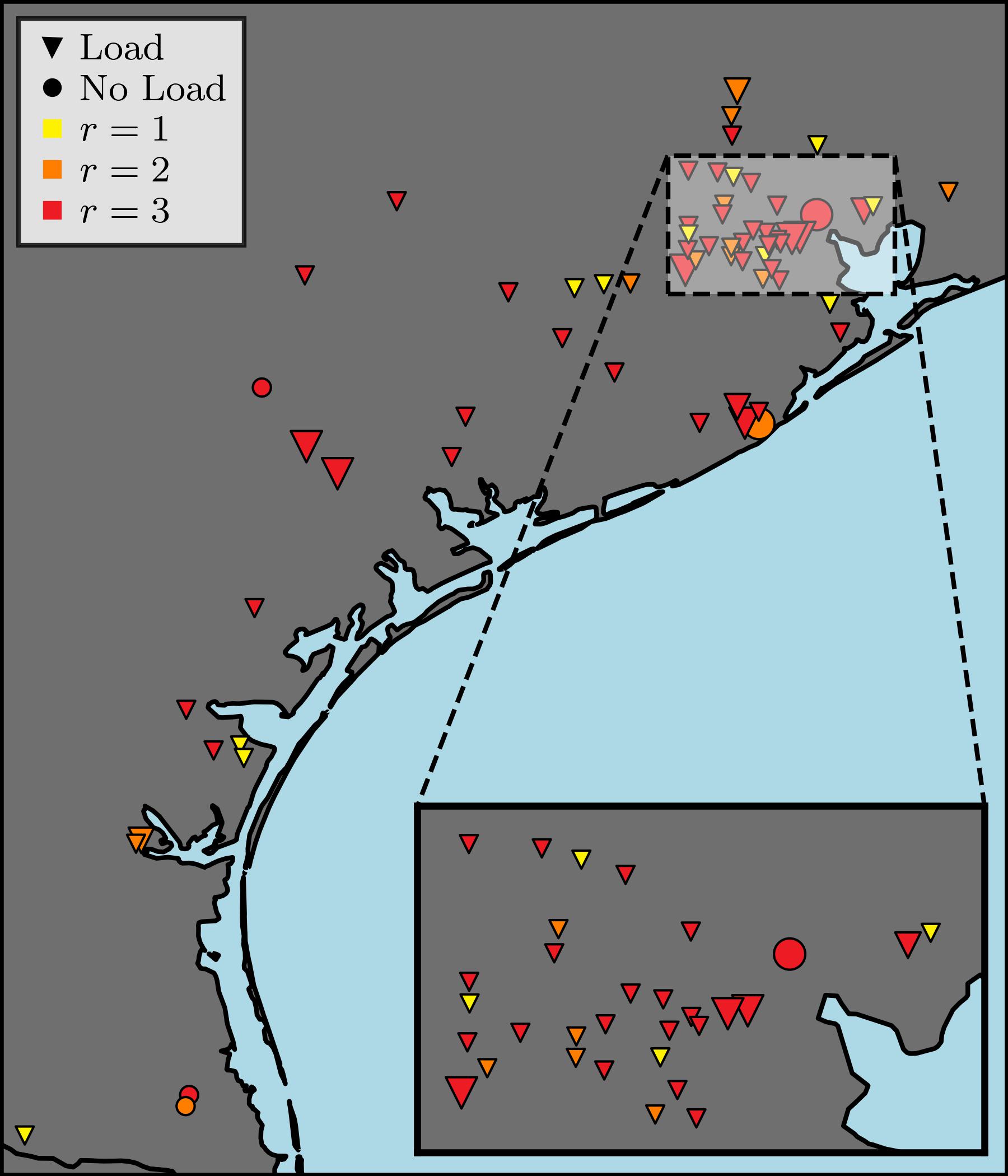}
        \caption{}
        \label{fig:harvey_uncertainty_map}
    \end{subfigure}
    \caption{In both the (a) Tropical Storm Imelda and (b) Hurricane Harvey scenarios, many substations experience inexorable flooding. In these maps, marker size corresponds to substation size, and marker color corresponds to the worst flood level experienced across all scenarios.}
\end{figure}

\begin{figure}
    \centering
    \includegraphics{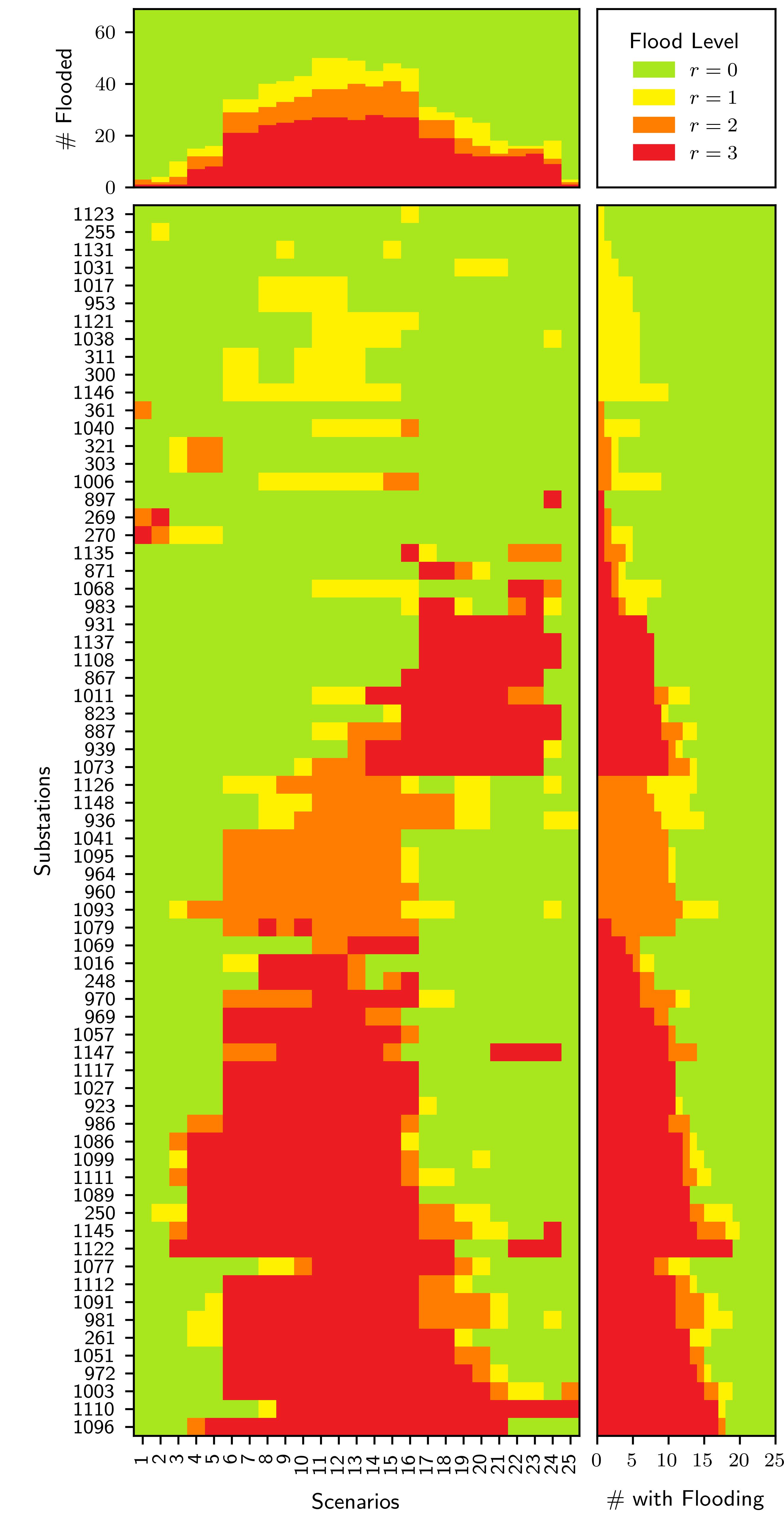}
    \caption{Hurricane Harvey scenario flood levels by scenario and by substation.}
    \label{fig:harvey_uncertainty_heatmap}
\end{figure}

The resulting flooding scenarios are illustrated in aggregate in Figure~\ref{fig:harvey_uncertainty_map} and an illustration of the full sample distribution in Figure~\ref{fig:harvey_uncertainty_heatmap}. Because the landfall locations are randomly sampled, we take the 25 scenarios to be equiprobable in our application.

%% file: results/_main.tex
\section{Results} \label{section:results}

\newcommand{\pathresults}{results}

\subsection{\emph{Experiments and Methodology}}

In this section, we assess the results collected from applying the two-stage models to the Tropical Storm Imelda and Hurricane Harvey case studies in a sensitivity study of the resource budget $f$. In both case studies, we set $\lambda^\text{shed} = \lambda^\text{over} = 1$. For the Imelda case study, we varied the budget between 0 and 20 in integer increments. For the Harvey case study, we studied integer budget values between 0 and 193. The precomputed upper limits of 20 and 193 are the maximum resources that may be used to effectively mitigate flooding.

Resilience metrics are generally either based on system performance (\textit{e.g.}, expected load shed) or system attributes (\textit{e.g.}, expected number of in-service components). Attribute-based metrics, though perhaps not as informative, are generally more simple to compute \citep{Vugrin2017}. To solve each instance, we first computed a handful of first-stage solutions using the parametric heuristic in Algorithm~\ref{alg:heuristics}. With generation capacity, transmission capacity, and serviceable load weighted by parameters $\eta^\text{gen}$, $\eta^\text{flow}$, and $\eta^\text{load}$, respectively, the algorithm greedily increases the expected weighted combination of these attributes. Because generation capacity is abundant in this application, we set $\eta^\text{gen} = 0$. We fixed $\eta^\text{load} = 1$ and varied ${\eta^\text{flow} \in \{0, 0.025, 0.05, 0.075, 0.1, 0.125, 0.15\}}$ to compute multiple solutions. The heuristic is iterative, and each iteration determines the allocation of resources to a singular substation that maximizes the ratio of expected marginal benefit to marginal cost. In the helper function $\textproc{Benefit}(\bs{x}, \tilde{\bs{x}}, \bs{\eta})$, $S(\cdot)$ embodies the logic of constraints \eqref{eq:dc_def_alpha} and \eqref{eq:dc_def_beta} that relate the mitigation decisions and flooding realizations to bus and branch statuses in each scenario. We implemented the extensive form of our model in Python using the \texttt{gurobipy} package and solved the instances using the Gurobi solver \citep{Gurobi2023} configured to leverage the heuristic solutions as partial warm-start solutions and then solve to optimality. We performed the optimization tasks on compute nodes with dual 48-core Intel\textregistered~Xeon\textregistered~Platinum 8160 CPUs from the Stampede2 cluster at the Texas Advanced Computing Center.
\begin{algorithm}
\caption{Parametric Greedy Heuristic} \label{alg:heuristics}
\begin{algorithmic}[1]
\Function{Heuristic}{$\bs{\eta}$, $f$}
    \State $\bs{x} \gets \bs{0}$
    \Comment $\bs{x} \in \{0,1\}^{|K \times R|}$
    \State $\displaystyle \mathcal{R} \gets \left\{\wt{\bs{x}} \in \mathcal{X}: \wt{\bs{x}} \ge \bs{x}, \exists! k \in K~\text{s.t.}~\sum_{r \in R} (\wt{x}_{kr} - x_{kr}) > 0\right\}$
    \While{$\displaystyle f > 0 \textbf{ and } \mathcal{R} \ne \emptyset \textbf{ and } \max_{\wt{\bs{x}} \in \mathcal{R}} \textproc{Benefit}(\bs{x}, \wt{\bs{x}}, \bs{\eta}) > 0$}
        \State $\displaystyle \wt{\bs{x}}^* \gets \argmax_{\wt{\bs{x}} \in \mathcal{R}} \textproc{Benefit}(\bs{x}, \wt{\bs{x}}, \bs{\eta}) / \bs{c}^\top (\wt{\bs{x}} - \bs{x})$
        \State $f \gets f - \bs{c}^\top (\wt{\bs{x}}^* - \bs{x})$
        \State $\bs{x} \gets \wt{\bs{x}}^*$
        \State $\displaystyle \mathcal{R} \gets \left\{\wt{\bs{x}} \in \mathcal{X}: \wt{\bs{x}} \ge \bs{x}, \exists! k \in K~\text{s.t.}~\sum_{r \in R} (\wt{x}_{kr} - x_{kr}) > 0\right\}$
    \EndWhile
    \State \Return $\bs{x}$
\EndFunction
\State
\Function{Benefit}{$\bs{x}$, $\wt{\bs{x}}$, $\bs{\eta}$}
    \State $\bs{\alpha}, \bs{\beta} \gets S(\bs{x}, \bs{\xi})$
    \Comment $\bs{\alpha} \in \{0,1\}^{|N \times \Omega|}, \bs{\beta} \in \{0,1\}^{|E \times \Omega|}$
    \State $\wt{\bs{\alpha}}, \wt{\bs{\beta}} \gets S(\wt{\bs{x}}, \bs{\xi})$
    \Comment $\wt{\bs{x}} \ge \bs{x} \implies \wt{\bs{\alpha}} \ge \bs{\alpha}, \wt{\bs{\beta}} \ge \bs{\beta}$
    \State $\displaystyle \rho^\text{load} \gets \sum_{\omega \in \Omega} \prob(\omega) \sum_{n \in N} (\wt{\alpha}_n^\omega - \alpha_n^\omega) \sum_{d \in D_n} p_d^\text{load}$
    \State $\displaystyle \rho^\text{gen} \gets \sum_{\omega \in \Omega} \prob(\omega) \sum_{n \in N} (\wt{\alpha}_n^\omega - \alpha_n^\omega) \sum_{g \in G_n} \overline{p}_g^\text{gen}$
    \State $\displaystyle \rho^\text{flow} \gets \sum_{\omega \in \Omega} \prob(\omega) \sum_{(n,m) \in E} (\wt{\beta}_{nm}^\omega - \beta_{nm}^\omega) \sum_{l \in L_{nm}} \overline{s}_l^\text{flow}$
    \State \Return $\rho^\text{load} \eta^\text{load} + \rho^\text{gen} \eta^\text{gen} + \rho^\text{flow} \eta^\text{flow}$
\EndFunction
\end{algorithmic}
\end{algorithm}

\begin{figure}
    \centering
    \ifdefined\epsswitch
        \includegraphics[width=6.5in]{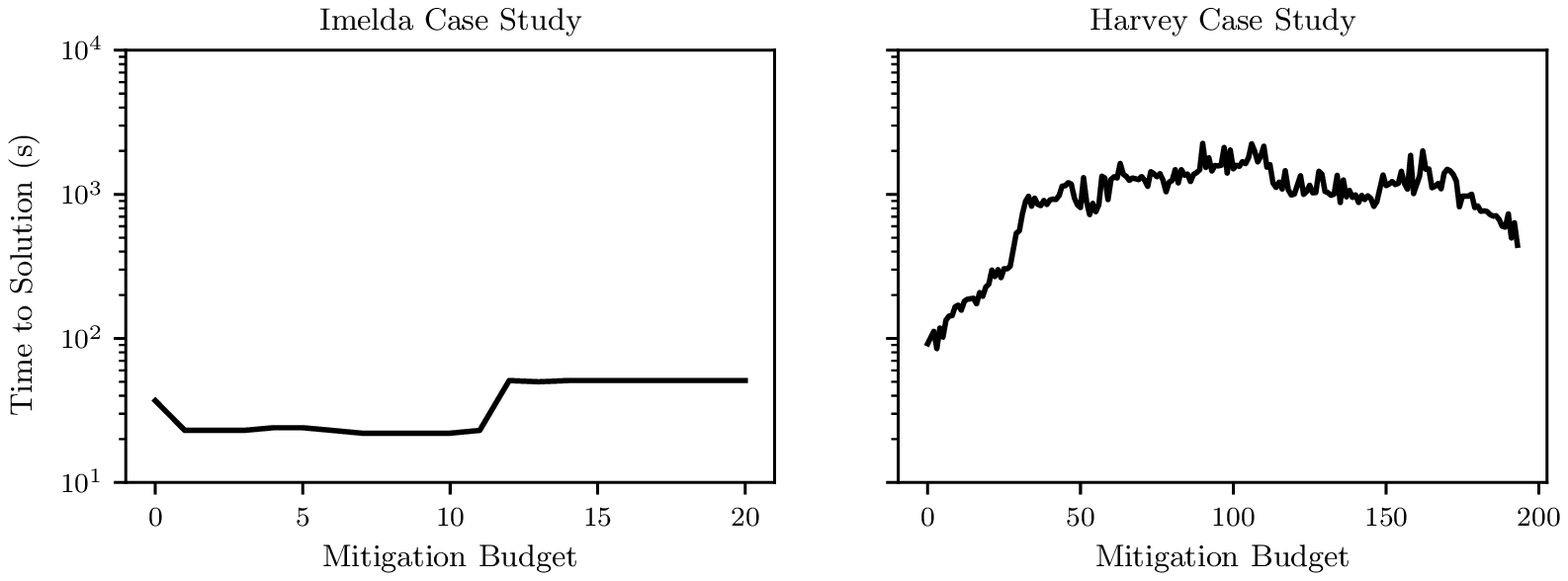}
    \else
        \includegraphics[width=6.5in]{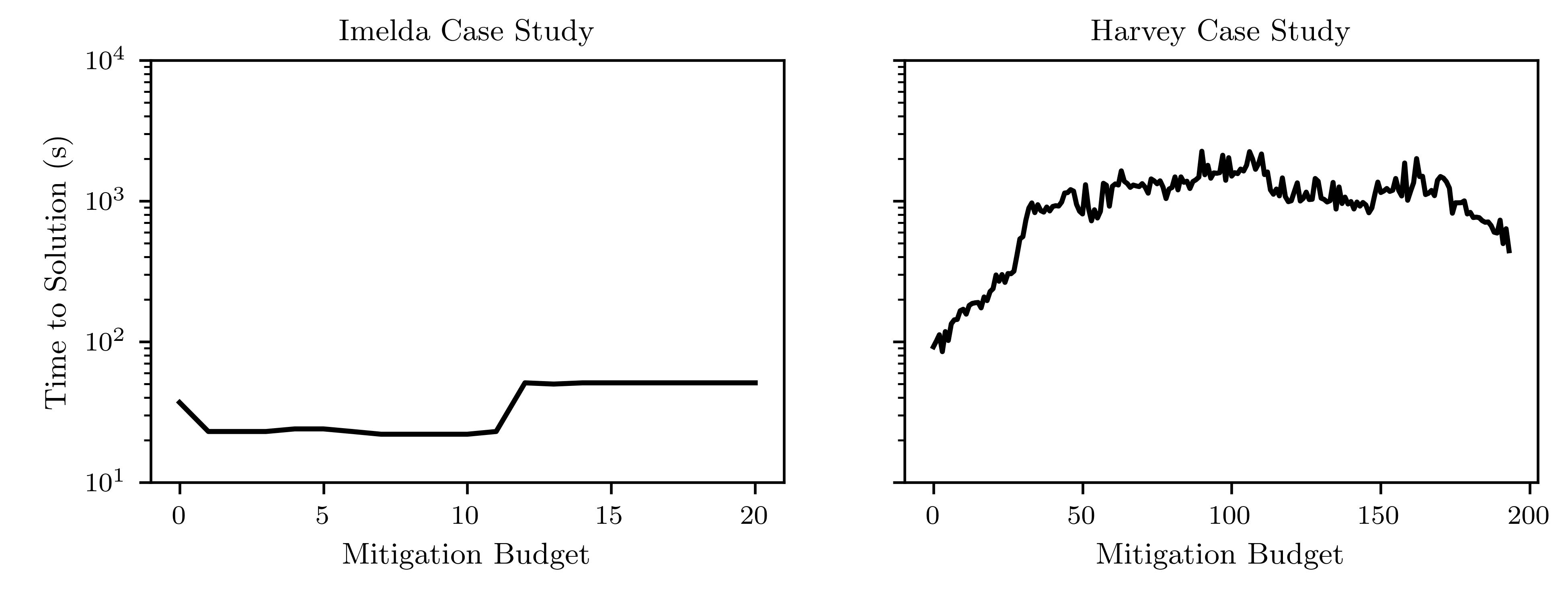}
    \fi
    \caption{Some of the Harvey instances required an order-of-magnitude more time to solve compared to the Imelda instances. These two plots share a logarithmically scaled vertical axis; however, the horizontal axis differs in limits and scale.}
    \label{fig:solution_times}
\end{figure}

The budget parameter $f$, though only directly affecting constraint \eqref{eq:con_resource_hi}, is crucial. The size of $\mathcal{X}$, the set of feasible first-stage solutions defined in \eqref{eq:resilience_X}, is strictly increasing in $f$. We use this insight to expedite the sensitivity study by solving the instances for each case study in increasing order of the budget and leveraging the obtained solutions as warm starts in subsequent instances. The solution times observed using this process to solve the SP model are shown in Figure~\ref{fig:solution_times}.

The relatively larger and more diverse set of scenarios from the Harvey case study account for the order-of-magnitude difference in the budget limits of 20 and 193 and similarly disparate solution times. Instances incorporating the Imelda scenario solve quickly, and the time-to-solution is insensitive to the resource budget. In contrast, instances incorporating the Harvey scenarios required between roughly 2 and 40 minutes to solve depending on the budget. Though the size of the first-stage decision space $\mathcal{X}$ is strictly increasing in the budget $f$ via constraint \eqref{eq:con_resource_hi}, we observe the longest times occurred for budgets of about 100 resources, roughly halfway between 0 and 193. We posit this occurs because the objective value is generally improved by deploying as many mitigation resources as possible, and the number of solutions using most if not all of the available resources begins to decrease as a function of the budget near $f = 100$. Ergo, the problem is most difficult to solve for intermediate budget values.

For the vast majority of the instances we studied, the optimal mitigation solutions were uniquely optimal. We determined this computationally by adding the no-good cut ${(\bs{1} - \bs{x}^*)^\top \bs{x} \ge 1}$ to the instance, solving that restricted instance to optimality, and assessing the objective value. In this cut, $\bs{x}^*$ is the identified optimal solution, $\bs{x}$ is the vector of decision variables, and $\boldsymbol{1}$ is the all-ones vector. For the vector of resource costs $\bs{c}$, if $\bs{c}^\top \bs{x}^* = f$, then the only solution cut from the feasible space is $\bs{x}^*$. Otherwise, this constraint cuts all $\bs{x} \ge \bs{x}^*$ (\textit{i.e.}, all the solutions that implement $\bs{x}^*$ and possibly more). Despite mitigating more flooding, such solutions are no better performing than $\bs{x}^*$ which, recall, is optimal.


We attribute the frequent occurrence of unique optimal mitigation solutions to the power grid instance having little symmetry. The few instances with multiple mitigation optima typically have one or more resources that cannot be used effectively and are instead deployed arbitrarily. As an example, consider a deterministic instance in which only one bus, a load bus, is affected. If 2 resources are required to prevent flooding at the associated substation, but only 1 resource is available, then that resource may be deployed arbitrarily with no effect. For the few instances having multiple mitigation optima, we did not bother identifying all optima as the task is computationally burdensome. That most instances have a unique mitigation optimum makes their comparison more meaningful.

\subsection{\emph{Analysis of Solutions}} \label{subsection:solution_analysis}
In the two-stage model, the mitigation decisions are the most important variables mainly because they must be made proactively. The other variables in the model, the power flow variables, represent decisions that are made reactively. Their inclusion in the model only serves to approximate the consequences of the mitigation in a set of representative scenarios. That said, we limit the scope of our analysis to the mitigation decisions.

\begin{figure}
    \centering
    \ifdefined\epsswitch
        \includegraphics[width=6.5in]{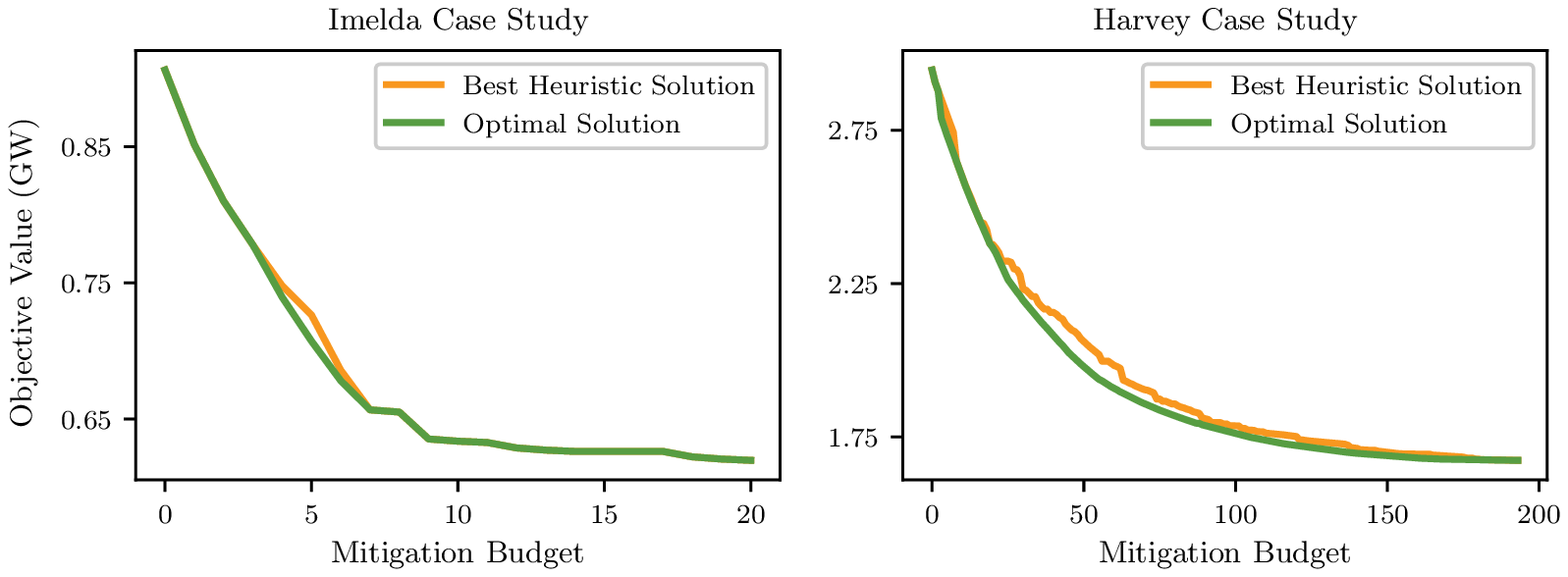}
    \else
        \includegraphics[width=6.5in]{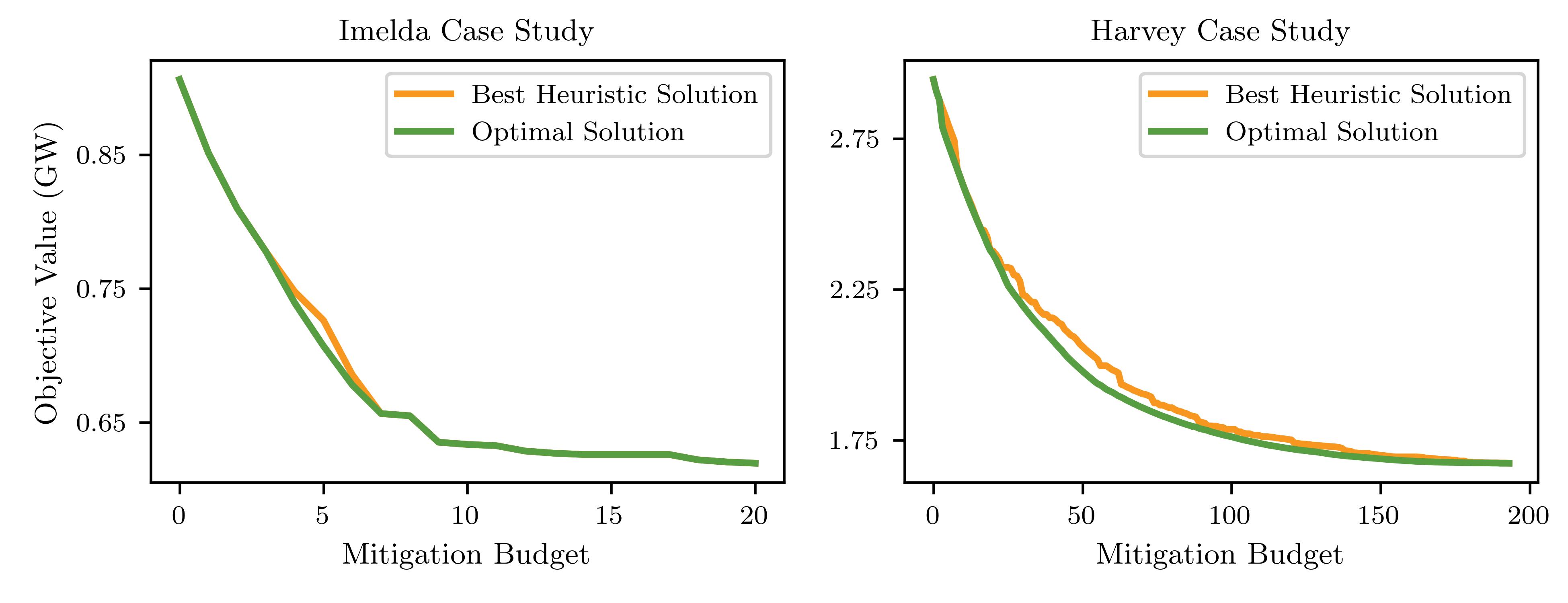}
    \fi
    \caption{In the Imelda and Harvey case studies, 31\% and 43\% decreases in the objective value are realized in the best objective value relative to the worst, and diminishing marginal returns are observed. Note that the axes in the two subplots are deliberately scaled differently to draw this comparison and that neither vertical axis begins at zero.}
    \label{fig:obj_and_bounds}
\end{figure}

Of course, the best measure of a mitigation solution's efficacy is its objective value. In our case studies, no overgeneration is necessary for the instances to be feasible, and our selected objective weights disincentivized it to the point that it was not present in any of our obtained optimal solutions. That is, the objective values represent load shed only.
\textcolor{response}{For each of the studied SP model instances, we present the objective values for the optimal solution and the best heuristic mitigation solution in Figure~\ref{fig:obj_and_bounds}.}
The two plots in the figure are scaled differently to draw a comparison. Specifically, we observe similarly diminishing marginal returns on the mitigation resources as the budget increases.
\textcolor{response}{For each instance, the relative optimality gap of the best heuristic mitigation solution $\bar{\boldsymbol{x}}$, measured by
\begin{equation*}
    \frac{
        \displaystyle \sum_{\omega \in \Omega} \text{Pr}(\omega) \mathcal{L}(\bar{\boldsymbol{x}}, \boldsymbol{\xi}^\omega) - \min_{\bs{x} \in \mathcal{X}} \sum_{\omega \in \Omega} \prob(\omega) \mathcal{L}(\bs{x}, \bs{\xi}^\omega)
    }{
        \displaystyle \min_{\bs{x} \in \mathcal{X}} \sum_{\omega \in \Omega} \prob(\omega) \mathcal{L}(\bs{x}, \bs{\xi}^\omega)
    },
\end{equation*}
is quite small. In fact, it is no more than 5\% for any of the instances.
Also,}
for any budget, the corresponding objective value ought to be considered relative to the worst-case and best-case consequences. The worst-case consequences occur when zero mitigation is enacted and the best-case when all preventable flooding is mitigated (\textit{i.e.}, the left- and right-most values, respectively, in each plot). In both case studies, notice that the best-case consequences are strictly positive due inexorable flooding sometimes affecting load buses. In the Imelda and Harvey case studies, 31\% and 43\% decreases in the objective value are realized in the best objective value relative to the worst.

As mentioned in the previous section, the size of the mitigation decision space $\mathcal{X}$ is non-decreasing in the mitigation budget $f$, and that is the only direct effect of the budget on the model. As such, the objective value is theoretically guaranteed to monotonically decrease, though perhaps not strictly, as the budget increases. For our two case studies, this is observed in Figure~\ref{fig:obj_and_bounds}.

We evaluated if the optimal mitigation solutions are similarly monotonic. Formally, we assessed if the optimal mitigation solutions satisfy $\bs{x}^*_f \le \bs{x}^*_{f+1}$ where $\bs{x}^*_f$ denotes an optimal solution for a budget $f$. In context, this equates to the optimal mitigation solutions being nested (\textit{i.e.}, optimal mitigation decisions for larger budgets are necessarily a superset of optimal decisions for smaller budgets with respect to both the location $k$ and resilience level $r$).

\begin{figure}
    \begin{subfigure}{\textwidth}
        \centering
        \includegraphics[width=6.5in]{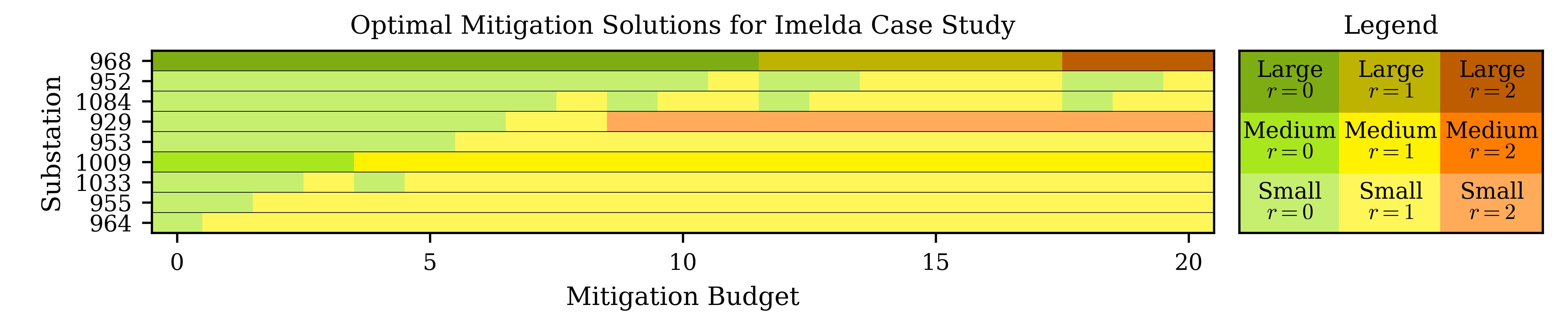}
        \caption{}
        \label{fig:solution-x-imelda}
    \end{subfigure}
    \par\bigskip
    \begin{subfigure}{\textwidth}
        \centering
        \includegraphics[width=6.5in]{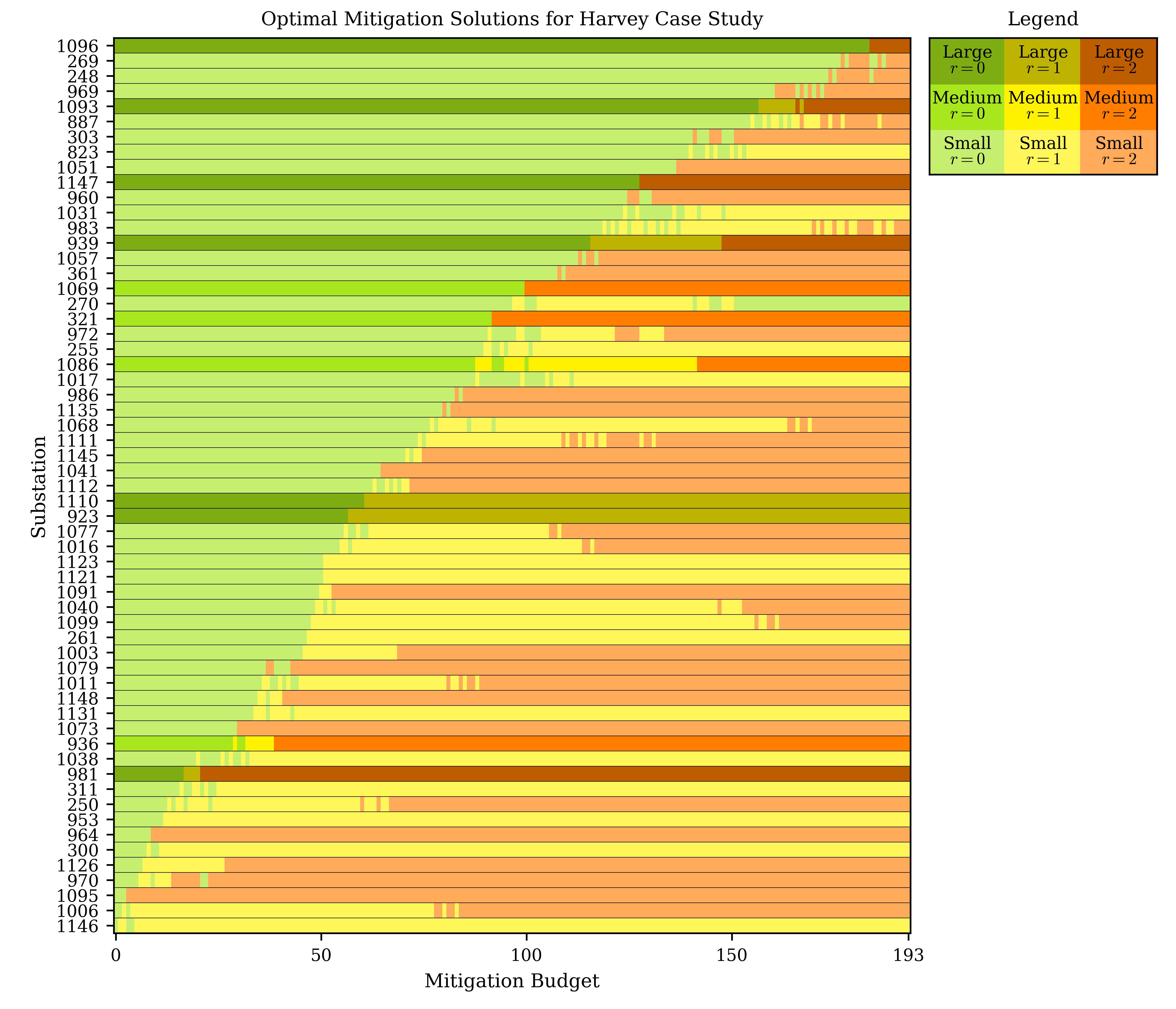}
        \caption{}
        \label{fig:solution-x-harvey}
    \end{subfigure}
    \caption{The mitigation solutions were not generally monotonic as a function of the budget in either (a) the Imelda case study or (b) the Harvey case study.}
    \label{fig:solution-x}
\end{figure}

The nature of discrete decision making is often such that optimal solutions are complex and even unintuitive. As an example, consider a toy 0/1 knapsack problem:
\begin{equation*}
    \max \left\{ 3 w_1 + 5 w_2 + w_3 : 4 w_1 + 8 w_2 + 3 w_3 \le C, \bs{w} \in \{0, 1\}^3 \right\}.
\end{equation*}
When $C = 7$, the optimal solution is $\bs{w}^* = [1, 0, 1]$. However, increment the budget by just one to $C = 8$ and $\bs{w}^* = [0, 1, 0]$ -- every decision flips! With consideration that $\mathcal{X}$ in our model features such complicating constraints, we did not expect the optimal solutions to be nested. Our identified optimal first-stage solutions, illustrated in Figure~\ref{fig:solution-x}, indeed did not exhibit nestedness. In those figures, color indicates the level of flooding to which a substation is made resilient through mitigation, and shading and tinting indicate substation size and thus the number of resources required for mitigation according to Table~\ref{tab:mitigation_costs}. Nonmonotonicity is visually evidenced by some rows comprising multiple continuous bands of the same color.

These figures highlight that the optimal mitigation for larger substations tends to be less sensitive to the budget than the optimal mitigation for smaller substations. In Figure~\ref{fig:solution-x-harvey}, for example, the optimal mitigation for large substation 939 transitions monotonically, but the mitigation for small substation 983 changes back-and-forth on 28 of the 193 studied budget increments. We also generally observe that the budget interval defined by a substation's first and last transitions from level $r$ to level $r + 1$ is short. That is, the sensitivity of the optimal mitigation solution to the budget at a specific substation is largely only transitory. However, because these transitory intervals are different at each substation, the solution exhibits nonmonotonicity over the entire range of studied budgets.

\begin{figure}[H]
    \centering
    \includegraphics[width=6.5in]{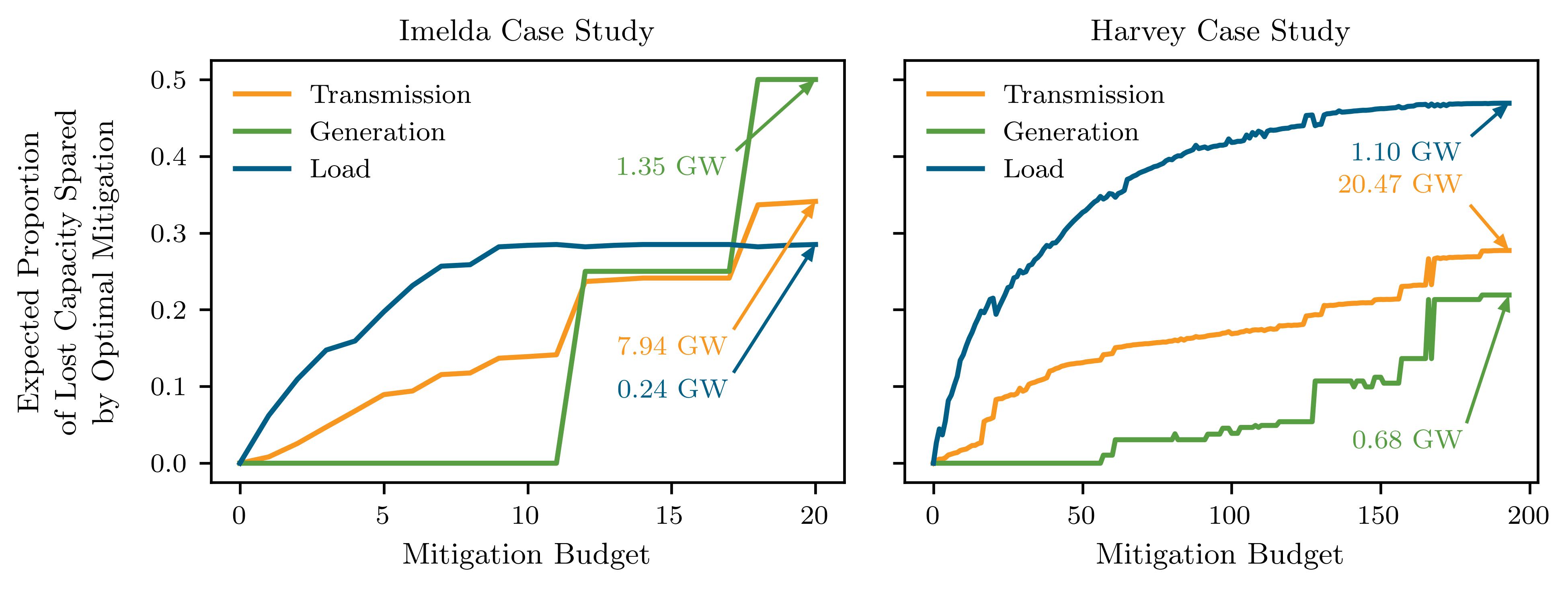}
    \caption{The expected load, generation, and transmission capacities spared by optimal mitigation are not generally monotonic. The illustrated quantities are as defined in \eqref{eq:expected_spared_load}, \eqref{eq:expected_spared_generation}, and \eqref{eq:expected_spared_transmission}. The annotations indicate the spared capacities in absolute units for the case in which all preventable flooding is mitigated.}
    \label{fig:improvement}
\end{figure}

Recall that equality constraints \eqref{eq:dc_def_alpha} and \eqref{eq:dc_def_beta} dictate the operational status of each component based on scenario flooding indicators $\bs{\xi}^\omega$ and enacted mitigation $\bs{x}$. As such, the nonmonotonicity of optimal mitigation propagates directly to load, generation, and transmission via $\bs{\alpha}^\omega$ and $\bs{\beta}^\omega$. We define ``lost capacity'' as that which is surrendered to flooding when no mitigation is enacted and ``spared capacity'' as that which would be lost if not for intervention. Normalized and expected measures of spared load, generation, and transmission capacities are
\begin{gather}
\sum_{\omega \in \Omega} \prob(\omega) \left(
    \bs{\alpha}^\omega - \accentset{\circ}{\bs{\alpha}}^\omega)^\top \bs{p}_n^\text{load} \middle/
    (\bs{1} - \accentset{\circ}{\bs{\alpha}}^\omega)^\top \bs{p}_n^\text{load}
\right), \label{eq:expected_spared_load} \\
\sum_{\omega \in \Omega} \prob(\omega) \left(
    \left(\bs{\alpha}^\omega - \accentset{\circ}{\bs{\alpha}}^\omega\right)^\top \overline{\bs{p}}_n^\text{gen} \middle/
    \left(\bs{1} - \accentset{\circ}{\bs{\alpha}}^\omega\right)^\top \overline{\bs{p}}_n^\text{gen}
\right), \label{eq:expected_spared_generation} \\
\sum_{\omega \in \Omega} \prob(\omega) \left(
    (\bs{\beta}^\omega - \accentset{\circ}{\bs{\beta}}^\omega)^\top \overline{\bs{s}}^\text{flow} \middle/
    (\bs{1} - \accentset{\circ}{\bs{\beta}}^\omega)^\top \overline{\bs{s}}^\text{flow}
\right). \label{eq:expected_spared_transmission}
\end{gather}
In these expressions, $\accentset{\circ}{\bs{\alpha}}^\omega$ and $\accentset{\circ}{\bs{\beta}}^\omega$ denote operational statuses that result in scenario $\omega$ if no mitigation is enacted, and $\bs{\alpha}^\omega$ and $\bs{\beta}^\omega$ denote those that result from enacting mitigation. In Figure~\ref{fig:improvement}, we refer to these expressions as the expected proportion of lost capacity spared by mitigation.
These figures highlight the nonmonotonicity of the expected spared capacities as functions of the budget. Note that the nonmonotonicity of the mitigation only guarantees that the spared capacities in each scenario are likewise nonmonotonic. In expectation, there is no such guarantee; however, our empirical results exhibit nonmonotonicity even in expectation.

A potentially problematic implication of the optimal mitigation being unnested is that the availability of more resources may lead to worse outcomes for some communities. This implication follows mainly from the effects of mitigation on the spared load and transmission capacities since there is an abundance of generation capacity in our two case studies. As an example, consider again substation 983 in the Harvey case study. For budget $f = 119$, the identified optimal solution involves implementing level $r = 1$ resilience at the substation. For budget $f = 120$, a larger budget, the identified optimal solution prescribes no mitigation there. Individuals whose load is served through that substation might perceive this as unfair and unintuitive. Of course, a net benefit is still realized by other communities having their loads spared and served instead.

\begin{figure}
    \begin{subfigure}{0.49\textwidth}
        \centering
        \includegraphics[width=3in]{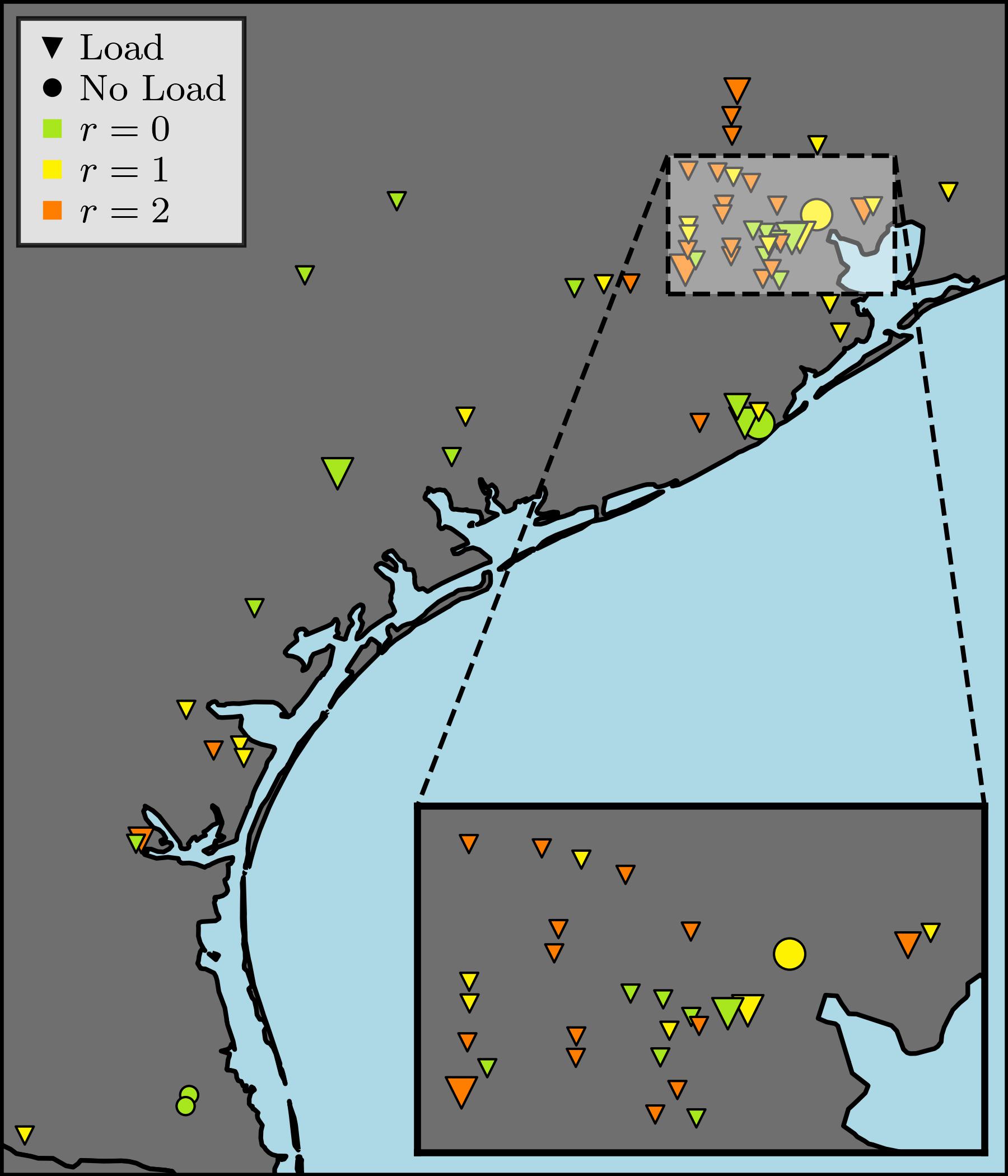}
        \caption{}
        \label{fig:harvey-map-budget-r3-f100}
    \end{subfigure}
    \begin{subfigure}{0.49\textwidth}
        \centering
        \includegraphics[width=3in]{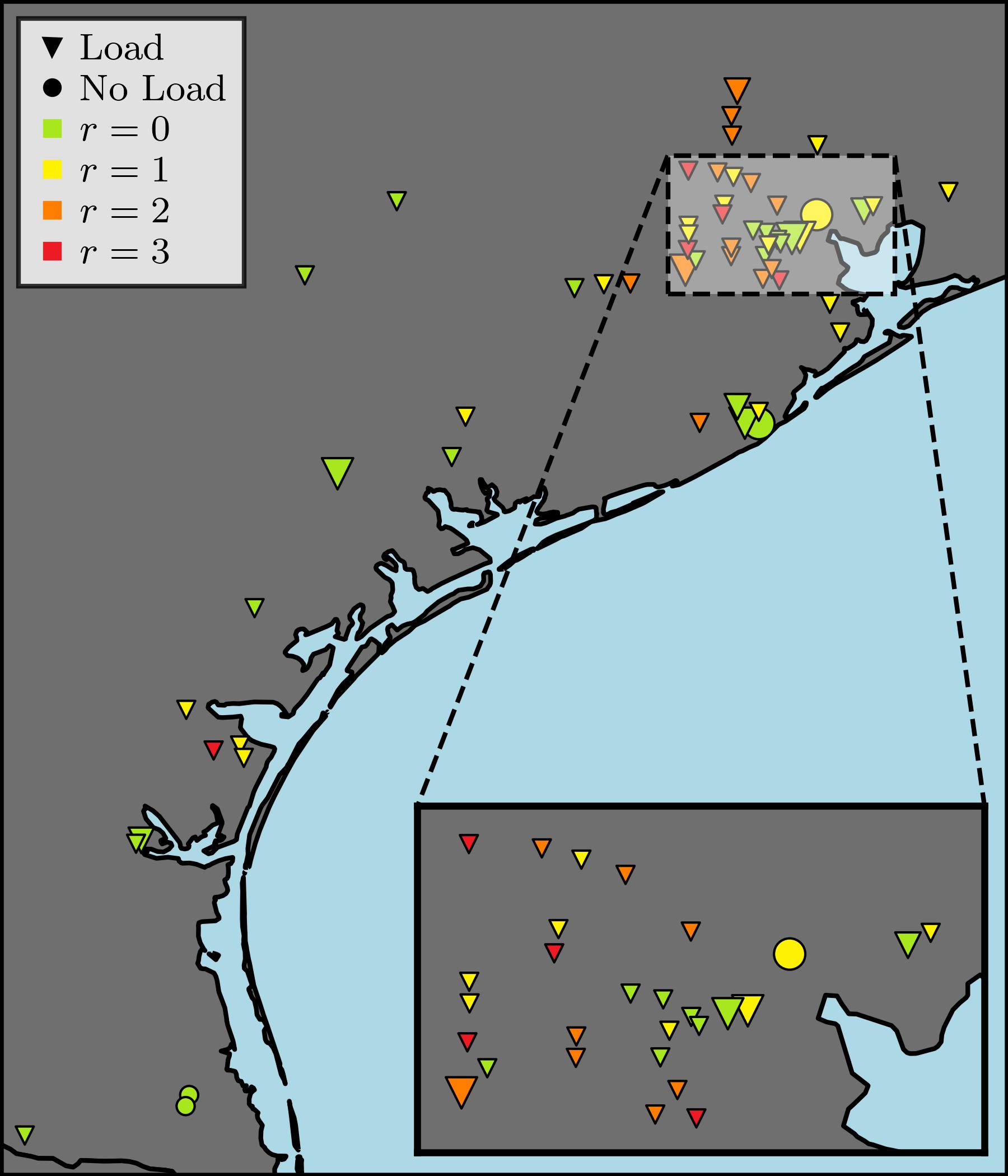}
        \caption{}
        \label{fig:harvey-map-budget-r4-f100}
    \end{subfigure}
    \caption{In the Harvey case study, the identified optimal mitigation differs considerably for budget $f = 100$ when the unattainable level of flooding is increased from (a) $\hat{r} = 3$ to (b) $\hat{r} = 4$. Marker size indicates substation size, and marker color indicates the optimal mitigation level.}
    \label{fig:rhat-sensitivity}
\end{figure}

Until now, all of the presented results have been based on instances with the unattainable level of flooding set to $\hat{r} = 3$. To conclude the results, we assess the impact of increasing that limit from $\hat{r} = 3$ to $\hat{r} = 4$ and discuss the spatial features of the optimal mitigation. We focus only on Harvey for this analysis since the corresponding set of scenarios is larger and more varied. In Figure~\ref{fig:rhat-sensitivity}, we contrast optimal mitigation solutions for for the cases of both $\hat{r} = 3$ and $\hat{r} = 4$ when the budget is fixed to $f = 100$.

Granted the budget-related transitory effects previously described are not captured by these illustrations, we see in both solutions that a large number of resources are allocated to substations in the more densely populated areas of Houston and Corpus Christi. Far fewer resources are allocated to the less densely populated areas elsewhere. This spatial disparity may be attributed to three main factors: flooding frequency and severity, load magnitudes, and grid physics and topology. It is immediately intuitive that substations with larger loads, worse flooding, and higher chances of flooding receive priority. Grid physics and topology affect decision making in varied ways, some more intuitive than others. One good example of an intuitive topological effect is at substation 1110, depicted by the large yellow circle to the northwest of Galveston Bay in Figure~\ref{fig:rhat-sensitivity}. This substation serves no load and suffers preventable flooding in only 1 scenario when $\hat{r} = 3$. Nevertheless, it is incident to 7 transmission lines, and its ability to facilitate the flow of power apparently warrants it receiving 3 resources for all budgets $f \ge 61$ as observed in Figure~\ref{fig:solution-x-harvey}.

Fixing $f$ and incrementing $\hat{r}$ causes the optimal mitigation to change in ways similar to how it changes when $\hat{r}$ is fixed and $f$ is incremented. In transition from Figure~\ref{fig:harvey-map-budget-r3-f100} to Figure~\ref{fig:harvey-map-budget-r4-f100}, resources are withdrawn from several substations so that other substations may be bolstered to level $r = 3$ for a net benefit: the optimal objective value improves from 1.76~GW to 1.66~GW. This is a roughly 8.5\% improvement relative to when no mitigation is deployed; however, the instance with $\hat{r} = 4$ requires about twice as much time to solve since the associated feasible mitigation solution set $\mathcal{X}$ is substantially larger. Because the marginal costs of mitigation levels are increasing, there is no difference in the optimal solution for sufficiently small budgets. However, the load shed reduced by unlimited mitigation resources improves from 43\% for $\hat{r} = 3$ with $f = 193$ to nearly 52\% for $\hat{r} = 4$ with $f = 304$.

%% file: conclusions/_main.tex
\section{Conclusions} \label{section:conclusions}

In this paper, we presented a two-stage stochastic programming model for informing power grid flood mitigation decision making prior to an imminent and uncertain hurricane. Our model of flood mitigation improved upon past models by capturing the flexibility and limitations of temporary flood barriers. To capture the consequences of flood barrier deployment decisions and unmitigated flooding on the power grid, we leveraged the DC power flow approximation in the second-stage recourse problems.
We applied these models to a pair of case studies featuring the ACTIVS 2000-bus synthetic grid of Texas and geographically realistic flooding scenarios derived from historical Tropical Storm Imelda and Hurricane Harvey data.

Our model is designed to inform how to deploy a fixed number of on-hand resources prior to a hurricane's imminent landfall. Assessing our results, we generally observed decreasing returns in the mitigation budget. This suggests that our model could be adapted to determine the number of resources that best balance the cost trade-offs of proactive mitigation and subsequent consequences supposing those costs are known.

We leveraged our results to highlight nonmonotonic trends in the optimal mitigation that arise as a result of the first-stage mitigation deployment problem comprising only discrete decision variables. In our model, increasing the budget broadens the set of feasible mitigation solutions. As such, the overall expected cost is guaranteed to decrease monotonically as a function of the budget. A potentially problematic implication, however, is that the net expected benefits achieved by deploying more resources to certain substations may come at the expense of substations that were protected at lower budget levels. Additionally, we highlight that the optimal mitigation solutions tend to allocate more resources to urban substations than to rural substations. Such solutions are preferred because they protect substations that experience worse or more likely flooding, serve larger loads, or otherwise facilitate the flow of power from generators to loads. Lastly, we show that permitting taller dams to be erected may lead to overall better outcomes but with optimal mitigation solutions that are quite different qualitatively. The benefit of allowing better mitigation in the model must be weighed against the detriment of the model being more difficult to solve.

In the future, we expect to improve the model in two ways. First, we would like to extend the recourse problem to include multiple time periods. Presently, we rely on performance under a single representative demand profile. Using a multi-time period model would allow us to more accurately capture temporal variability like the daily periodicity of demand and renewable generation and temporal constraints like generator ramp rate limits. We believe this would improve the applicability of the prescribed solutions to reality. Of course, this would increase the complexity of the recourse problem and assuredly require more time to solve the model. Should it prove to be too difficult a model to solve, we could at the very least attempt to validate the prescribed solutions in multi-time period scenarios. Second, resilience in this paper was quantified by expected system-wide load shed, but other objectives or side constraints related to performance may be considered. Alternatives such as the joint minimization of mitigation and load shed costs, equity-related metrics, service level constraints, chance constraints, etc. may be considered in the future. \textcolor{response}{Regarding analysis, assessing the sensitivity of our models to the uncertainty parameterization would also be prudent. Understanding how the prescribed mitigation is affected by varying the scenario sampling scheme (\textit{e.g.}, by including or omitting stratification, adopting importance sampling, adjusting the number of sampled scenarios, and moreover evaluating if the most effective scheme depends on other problem data like the mitigation limits and power grid parameters) would help inform how the models may be most effectively used by practitioners.}